\documentclass[numreferences]{kluwer}
\usepackage[all]{xy}
\newtheorem{thm}{Theorem}
\newtheorem{cor}[thm]{Corollary}
\newtheorem{prop}[thm]{Proposition}
\newtheorem{lem}[thm]{Lemma}
\newdisplay{defin}[thm]{Definition}
\newdisplay{rmk}[thm]{Remark}
\newdisplay{exam}[thm]{Example}
\newdisplay{case}[thm]{Case}

\newcommand{\mb}[1]{\mathbb{#1}}

\newcommand{\Hom}{\ensuremath{{\rm Hom}}}
\newcommand{\colim}{\ensuremath{\mathop{\rm colim}}}

\newcommand{\hocolim}{\ensuremath{\mathop{\rm hocolim}}}

\newcommand{\QIrr}{\ensuremath{\mathop{\rm QIrr}}}

\newcommand{\Ob}{\ensuremath{\mathop{\rm Ob}}}
\newcommand{\Mor}{\ensuremath{\mathop{\rm Mor}}}
\newcommand{\Sing}{\ensuremath{\mathop{\rm Sing}}}
\newcommand{\U}{{\rm U}}
\newcommand{\PU}{{\rm PU}}
\newcommand{\BU}{{\rm BU}}
\newcommand{\EU}{{\rm EU}}
\newcommand{\Rig}{{\rm Rig}}

\newcommand{\eilm}[1]{\ensuremath{{\rm H} #1}}
\newcommand{\smsh}[1]{\ensuremath{\mathop{\wedge}_{#1}}}

\newcommand{\pull}[1]{\ensuremath{\mathop{\times}_{#1}}}

\newcommand{\cat}[1]{\ensuremath{{\bf #1}}}
\newcommand{\multc}[1]{\ensuremath{{\bf #1}}}

\newcommand{\Gm}[2]{\ensuremath{\Gamma^o_G(#1,#2)}}

\newcommand{\xym}[1]{
\vskip 0.7pc
\centerline{\xymatrix{
#1
}}
\vskip 0.7pc
\noindent}

\begin{document}

\begin{article}
\begin{opening}
\title{The product formula in unitary deformation $K$-theory}
\author{Tyler
  \surname{Lawson}\email{tlawson@math.mit.edu}\thanks{Partially
    supported by NSF grant 0402950.}}
%\runningauthor{Tyler Lawson}
%\runningtitle{The product formula in unitary deformation $K$-theory}
\institute{Massachusetts Institute of Technology, Cambridge, MA 02139}

\begin{abstract}
  For finitely generated groups ${\rm G}$ and ${\rm H}$, we prove
  that there is a weak equivalence ${\cal K}{\rm G} \smsh{ku} {\cal K}{\rm H}
  \simeq {\cal K}({\rm G} \times {\rm H})$ of $ku$-algebra spectra,
  where ${\cal K}$ denotes the ``unitary deformation $K$-theory''
  functor.  Additionally, we give spectral sequences for computing the
  homotopy groups of ${\cal K}{\rm G}$ and $\eilm{\mb Z} \smsh{ku} {\cal
  K}{\rm G}$ in terms of connective $K$-theory and homology of
  spaces of ${\rm G}$-representations.
\end{abstract}
%\keywords{}
\end{opening}

\section{Introduction}

The underlying goal of many programs in algebraic $K$-theory is to
understand the algebraic $K$-groups of a field $F$ as being built from
the $K$-groups of the algebraic closure of the field, together with
the action of the absolute Galois group.  Specifically, Carlsson's
program (see~\cite{carlsson:kth}) is to construct a model for the
algebraic $K$-theory {\em spectrum\/} using the Galois group and the
$K$-theory spectrum of the algebraic closure $\overline F$.

In some specific instances, the absolute Galois group of the field $F$
is explicitly the profinite completion $\hat {\rm G}$ of a discrete
group ${\rm G}$.  (For example, the absolute Galois group of the
field $k(z)$ of rational functions, where $k$ is an algebraically
closed of characteristic zero, is the profinite completion of a free
group.)  In the case where $F$ contains an algebraically closed
subfield, the profinite completion of a ``deformation $K$-theory''
spectrum ${\cal K}{\rm G}$ is conjecturally equivalent to the profinite
completion of the algebraic $K$-theory spectrum ${\mb K}F$.

Additionally, it would be advantageous for this description to be
compatible with the motivic spectral sequence.  This deformation
$K$-theory spectrum has an Atiyah-Hirzebruch spectral sequence arising
from a spectrum level filtration.  The filtration quotients are
spectra built from isomorphism classes of representations of the
group.  It is hoped that this filtration is related to the motivic
spectral sequence, and that this relation would give a greater
understanding of the relationships between Milnor $K$-theory, Galois
cohomology, and the representation theory of the Galois group.

We now outline the construction of deformation $K$-theory.  To a finitely
generated group ${\rm G}$ one associates the category ${\cal C}$ of
finite dimensional unitary representations of ${\rm G}$, with
morphisms being equivariant isometric isomorphisms.  Elementary
methods of representation theory allow this category to be analyzed;
explicitly, the category of unitary ${\rm G}$-representations is
naturally equivalent to a direct sum of copies of the (topological)
category of unitary vector spaces.

However, there is more structure to ${\cal C}$.  First, there is a
bilinear tensor product pairing.  Second, ${\cal C}$ can be given the
structure of an internal category ${\cal C}^{top}$ in $\cat{Top}$.
This means that there are {\em spaces\/} $\Ob({\cal C}^{top})$ and
$\Mor({\cal C}^{top})$, together with continuous domain, range,
identity, and composition maps, satisfying appropriate associativity
and unity diagrams.  The topology on $\Ob({\cal C}^{top})$ reflects
the possibility that homomorphisms from ${\rm G}$ to $\U(n)$ can
continuously vary from one isomorphism class of representations to
another.  The identity map on objects and morphisms is a continuous
functor ${\cal C} \to {\cal C}^{top}$ that is bijective on objects.

Both of these categories have notions of direct sums and so are
suitable for application of an appropriate infinite loop space
machine, such as Segal's machine~\cite{catcohom}.  This yields a map
of (ring) spectra as follows:
\begin{equation}
  \label{eq:ringmap}
\mb K {\cal C} \simeq \bigvee ku \to \mb K{\cal
  C}^{top} \ \mathop{=}^{\rm def}\  {\cal K}{\rm G}.
\end{equation}
Here $\mb K$ is an algebraic $K$-theory functor, $ku$ is the
connective $K$-theory spectrum, and the wedge is taken over the set of
irreducible unitary representations of ${\rm G}$.  Note that
$\pi_*(\mb K{\cal C}) \cong R[{\rm G}] \otimes \pi_*(ku)$ as a
ring, where $R[{\rm G}]$ is the unitary representation ring of ${\rm
G}$.

The spectrum ${\cal K}{\rm G}$ is the {\em unitary deformation $K$-theory\/}
of ${\rm G}$.  It differs from the $C^*$-algebra $K$-theory of ${\rm
G}$---for example, in section~\ref{sec:exactcouple}, we find that
the unitary deformation $K$-theory of the discrete Heisenberg group
has infinitely generated $\pi_0$.

When ${\rm G}$ is free on $k$ generators, one can directly verify the
formula
\[
{\cal K}{\rm G} \simeq ku \vee \left(\bigvee^k \Sigma ku\right).
\]
A more functorial description in this case is that ${\cal K}{\rm G}$ is the
connective cover of the function spectrum $F(B{\rm G}_+,ku)$.  This
formula does not hold in general; even in simple cases ${\cal K}{\rm G}$ can
be difficult to directly compute, such as when ${\rm G}$ is free
abelian on multiple generators.

In this paper, we will prove the following product formula for unitary
deformation $K$-theory.

\begin{thm}
\label{thm:product}
The tensor product map induces a map of commutative $ku$-algebra
spectra ${\cal K}{\rm G} \smsh{ku} {\cal K}{\rm H} \to {\cal K}({\rm G} \times
{\rm H})$, and this map is a weak equivalence.
\end{thm}

The reader should compare the formula
\[
R[{\rm G}] \otimes R[{\rm H}] \cong R[{\rm G} \times {\rm H}]
\]
for unitary representation rings.

The proof of Theorem~\ref{thm:product} proceeds by making use of a
natural filtration of ${\cal K}{\rm G}$ by subspectra ${\cal K}{\rm G}_n$.  These
subspectra correspond to representations of ${\rm G}$ whose
irreducible components have dimension less than or equal to $n$.
Specifically, we show in sections~\ref{sec:local} and
\ref{sec:identify} that there is a homotopy fibration sequence of
spectra
\[ {\cal K}{\rm G}_{n-1} \to {\cal K}{\rm G}_n \to 
\left(\Hom({\rm G},\U(n))/{\rm Sum}({\rm G},n)\right)
\smsh{\PU(n)} ku^{\PU(n)}.
\]
Here ${\rm Sum}({\rm G},n)$ is the subspace of $\Hom({\rm G},\U(n))$
consisting of those representations ${\rm G} \to \U(n)$ containing a
nontrivial invariant subspace.  The spectrum $ku^{\PU(n)}$ is a
connective $\PU(n)$-equivariant $K$-homology spectrum
discussed in section~\ref{sec:ctkh}.

As side benefits of the existence of this filtration,
Theorems~\ref{thm:kspecseq} and~\ref{thm:homspecseq} give spectral
sequences for computing the homotopy groups of ${\cal K}{\rm G}$ and the
homotopy groups of $\eilm{\mb Z} \smsh{ku} {\cal K}{\rm G}$ respectively.

When ${\rm G}$ is free on $k$ generators,
Theorem~\ref{thm:homspecseq} gives a spectral sequence converging to
$\mb Z$ in dimension 0, $\mb Z^k$ in dimension 1, and 0 otherwise, but
the terms in the spectral sequence are highly nontrivial---they are
the homology groups of the spaces of $k$-tuples of elements of
$\U(n)$, mod conjugation and relative to the subspace of $k$-tuples
that contain a nontrivial invariant subspace.  The method by which the
terms in this spectral sequence are eliminated is a bit mysterious.

The layout of this paper is as follows.  Section~\ref{sec:gamma} gives
the necessary background on $G$-equivariant $\Gamma$-spaces for a
compact Lie group $G$, allowing identification of equivariant smash
products.  The model theory of such functors was considered when $G$
is a finite group in~\cite{dro}, using simplicial spaces.  Our
approach to the proofs of the results we need follows the approach
of~\cite{dgm}.  In Section~\ref{sec:ctkh} we give explicit
constructions of an equivariant version of connective $K$-theory.  In
Section~\ref{sec:defk} the unitary deformation $K$-theory of ${\rm
G}$ is defined, and Sections~\ref{sec:local} and~\ref{sec:identify}
are devoted to constructing the localization sequences, in particular
explicitly identifying the base as an equivariant smash product.  In
section~\ref{sec:einfty} the algebra and module structures are made
explicit by making use of results of Elmendorf and Mandell~\cite{rma}.
The proofs of the main theorems are completed in
sections~\ref{sec:exactcouple} and~\ref{sec:product}.

A proof of the product formula for representations in ${\rm GL}(n)$,
rather than $\U(n)$, would also be desirable.  This paper makes use of
quite rigid constructions that make apparent the identification of the
base in the localization sequence with a particular model for the
equivariant smash product.  In the case of ${\rm GL}(n)$, the
definitions of both the cofiber in the localization sequence and the
equivariant smash product need to be replaced by notions that are more
well-behaved from the point of view of homotopy theory.

\section{Preliminaries on $G$-equivariant $\Gamma$-spaces}
\label{sec:gamma}

In this section, $G$ is a compact Lie group, and actions of $G$ on
based spaces are assumed to fix the basepoint; a {\em free} action
will be one that is free away from the basepoint.  We will now carry
out constructions of $\Gamma$-spaces in a na\"ive equivariant context.
When $G$ is trivial these are definitions for topological
$\Gamma$-spaces, as in the appendix of~\cite{schwede}.

For any natural number $k$, denote the based space
$\{*,1,\ldots,k\}$ by $k_+$.

Let $\Gamma^o_G$ be the category of {\it right} $G$-spaces that are
isomorphic to ones the form $G_+ \smsh{} k_+$, with morphisms being
$G$-equivariant.  (Strictly speaking, we take a small skeleton for
this category.)  The set $\Gm{X}{Y}$ can be given the mapping space
topology, giving this category an enrichment in spaces.  Explicitly,
\[
\Gm{G_+ \smsh{} k_+}{Z} \cong \prod^k Z
\]
as a space.  If $G$ is trivial we drop it from the notation.

\begin{defin}
A $\Gamma_G$-space is a base-point preserving continuous functor
$\Gamma^o_G \to \cat{Top}_*$.
\end{defin}

Here $\cat{Top}_*$ is the category spaces, i.e. compactly generated
weak Hausdorff pointed topological spaces with nondegenerate
basepoint, which has internal function objects $F(-,-)$.

Any $\Gamma_G$-space $M$ has an underlying $\Gamma$-space $M(G_+
\smsh{} -)$.  This $\Gamma$-space inherits a continuous left
$G$-action because the left action of $G$ on the first factor of $G_+
\smsh{} Y$ is right $G$-equivariant.  This can be expressed as a
continuous map
\[
\phi : G \to \Gm{G_+ \smsh{} Y}{G_+ \smsh{} Y},
\]
given by $\phi(g)(h \smsh{} y) = gh \smsh{} y$.  Then the composite
map
\[
M \circ \phi: G \to F(M(G_+ \smsh{} Y), M(G_+ \smsh{} Y)),
\]
gives an action of $G$ on $M(G_+ \smsh{} Y)$ that is natural in $Y$.
This gives a $G$-action on the underlying $\Gamma$-space of $M$.

\begin{rmk}
More generally, if $H \to G$ is a map of groups, the formula $M
\mapsto M(- \smsh{H} G_+)$ defines a restriction map from
$\Gamma_G$-spaces to $\Gamma_H$-spaces with a left action of $C(H)$,
the centralizer of $H$ in $G$.
\end{rmk}

For technical reasons, we require the following definition.

\begin{defin}
A $\Gamma_G$-space $M$ is semi-cofibrant if for any $Z \in
\Gamma^o_G$, the inclusion
\[
\tilde M (Z) = \bigcup_{Y \subsetneq Z} M(Y) \subset M(Z)
\]
is a cofibration of spaces.
\end{defin}

A cofibrant $\Gamma$-space as defined in~\cite{bousfieldfriedlander}
would require some freeness assumption about the action of ${\rm
Aut}(Z)$ on $M(Z)/\tilde M(Z)$.

\begin{rmk}
\label{rmk:degen}
For any element $m \in M(Z)$, there is a unique minimal subobject $Y$
of $Z$ such that $m$ is in the image of $M(Y)$.  To see this, consider
the inclusion $i_Y : Y \to Z$ and the map $\pi_Y : Z \to Y$ that is
the identity on $Y$ and sends the rest of $Z$ to the basepoint.
The map $i_Y \pi_Y$ acts as the identity on the image of $M(Y)$, and
$(i_Y \pi_Y) (i_{Y'} \pi_{Y'}) = i_{Y \cap Y'} \pi_{Y \cap Y'}$.
\end{rmk}

Let $X,Z \in \Gamma^o_G$, $Y$ a based set.  The space $G_+ \smsh{} Y$
has commuting left and right $G$-actions.  This gives rise to a
continuous right action of $G$ on $\Gm{G_+ \smsh{} Y}{Z}$.  There is a
map $\phi: X \to \Gm{G_+ \smsh{} Y}{X \smsh{} Y}$ given by $x \mapsto
\phi_x$, where $\phi_x(g \smsh{} y) = xg \smsh{} y$.  The map $\phi$
is clearly right $G$-equivariant.  Composing the map $\phi$ with the
functor $M$ gives a continuous $G$-equivariant map from $X$ to
$F(M(G_+ \smsh{} Y), M(X \smsh{} Y))$, and the adjoint is a natural
assembly map $X \smsh{G} M(G_+ \smsh{} Y) \to M(X \smsh{} Y)$.

We can promote a $\Gamma_G$-space $M$ to a functor on all free right
$G$-spaces, as follows.
\begin{defin}
For $X$ a right $G$-space, the $\Gamma$-space $X \otimes_G M$ is defined
by:
\[
X \otimes_G M(Z) = \left(\coprod_{Y \in \Gamma^o_G} 
M(Y) \smsh{} F^G(Y,X \smsh{} Z)
  \right) / \sim.
\]
Here the equivalence relation $\sim$ is generated by relations
$(u \smsh{} f^* v) \sim (f_* u \smsh{} v)$ for $f : Y \to Y', u \in
M(Y), v \in F^G(Y',X \smsh{} Z)$.  More concisely, $X \otimes_G
M(Z)$ can be expressed as the coend
\[
\int^Y M(Y) \smsh{} F^G(Y,X\smsh{}Z).
\]
\end{defin}

\begin{rmk}
  If $X_.$ is a simplicial object in the category $\Gamma^o_G$, there
  is a natural homeomorphism $|X_.| \otimes_G M \to |M(X_.)|$.  (A short
  proof can be given by expressing $|X_.|$ as a coend $\int^n X_n
  \smsh{} \Delta^n_+$ and applying ``Fubini's
  theorem''---see~\cite{cwm}, chapter IX.) The reason for allowing
  $G$-CW complexes rather than simply restricting to these simplicial
  objects is that some $G$-homotopy types cannot be realized by
  simplicial objects.  For example, any such simplicial object of the
  form $X_+$ is a principal $G$-bundle over $X/G_+$, and is classified
  by an element in $H^1_{discrete}(X/G,G)$.  A general $X_+$ is
  classified by an element in $H^1_{cont.}(X/G,G)$.
\end{rmk}

We will refer to the space $(X \otimes_G M)(1_+)$ as $M(X)$; this
agrees with the notation already defined when $X \in \Gamma^o_G$.  We
will only apply this construction to cofibrant objects in a certain
model category of based $G$-spaces; specifically, we will only apply
this construction to based $G$-CW complexes with free action away from
the basepoint.  Such objects are formed by iterated cell attachment of
$G_+ \smsh{} D^n_+ $ along $G_+ \smsh{} S^{n-1}_+$.

It will be useful to have homotopy theoretic control on $X \otimes_G
M$, for the purposes of which we introduce a less rigid tensor product.

\begin{defin}
For $M$ a $\Gamma_G$-space, we can define a simplicial
$\Gamma_G$-space $LM_.$ by setting $LM(Z)_p$ equal to
\[
\bigvee_{Z_0,\ldots,Z_p \in \Gamma^o_G} \Gm{Z_p}{Z} \smsh{}
\Gm{Z_{p-1}}{Z_p} \smsh{} \cdots \smsh{} \Gm{Z_0}{Z_1}
\smsh{} M(Z_0).
\]
The face maps are given by:
\begin{eqnarray*}
  d_i(f_0\smsh{}\ldots\smsh{}f_p\smsh{}m) &=& f_0 \smsh{} \ldots
  \smsh{} f_i \circ f_{i+1} \smsh{} \ldots \smsh{} f_p \smsh{}
  m\hbox{\ \ if $i<p$,}\\
  d_p(f_0\smsh{}\ldots\smsh{}f_p\smsh{}m) &=& f_0 \smsh{} \ldots
  \smsh{} f_{p-1} \smsh{} (Mf_p)(m).
\end{eqnarray*}
The degeneracy map $s_i$ is an insertion of an identity map after
$f_i$ for $0 \leq i \leq p$.
\end{defin}

Suppose $M$ is a semi-cofibrant $\Gamma_G$-space.  For any $p$ and $Z$, the
subspace of $LM(Z)_p$ consisting of degenerate objects is the union of
the subspaces that contain an identity element in some component of
the smash product. All the spaces in the smash product are
nondegenerately based cofibrant objects, and so the inclusion of the
degenerate subcomplex is a cofibration.  As a result, this simplicial
space is {\em good\/} in the sense of~\cite{catcohom}, Appendix A, and
so the geometric realization of it is homotopically well-behaved.

The simplicial $\Gamma_G$-space $LM_.$ has a natural augmentation $LM_. \to
M$.  The augmented object $LM_. \to M$ has an extra degeneracy map
$s_{-1}$, defined by
\[
s_{-1}(f_0 \smsh{} \ldots \smsh{} f_p \smsh{} m) =
id \smsh{} f_0 \smsh{} \ldots \smsh{} f_p \smsh{} m.
\]
As a result, the map $|LM_.(Z)| \to M(Z)$ is a homotopy equivalence for
any $Z \in \Gamma^o_G$.  (Note that $LM_.$ is the bar construction
$B(\Gamma^o_G, \Gamma^o_G, M).$)

For $X$ a right $G$-space, consider the simplicial $\Gamma_G$-space $X
\otimes_G LM_.$.  We have
\[
X \otimes_G LM_p = \bigvee_{Z_0,\ldots,Z_p} \Big[X \otimes_G
\Gm{Z_p}{-}\Big] \smsh{} \Gm{Z_{p-1}}{Z_p} \smsh{} \cdots
\smsh{} M(Z_0)
\]
as a $\Gamma$-space, because the tensor construction distributes over
wedge products and commutes with smashing with spaces.  However, a
straightforward calculation yields the formula
\[
\Big[X \otimes_G \Gm{Y}{-}\Big](Z) \cong F^G(Y,X \smsh{} Z),
\]
the space of $G$-equivariant based functions from $Y$ to $X \smsh{}
Z$.

\begin{prop} For $M$ a semi-cofibrant $\Gamma_G$-space and $X$ a $G$-CW
  complex with free action away from the basepoint, the augmentation
  map $X \otimes_G LM_. \to X \otimes_G M$ is a levelwise weak
  equivalence of $\Gamma$-spaces after realization; i.e., the map
  $|(X \otimes_G LM_.)(Z)| \to (X \otimes_G M)(Z)$ is a weak equivalence
  for all $Z \in \Gamma^o$.
\end{prop}

\begin{pf}
There is a natural isomorphism $(X \otimes_G M)(Z) \cong M(X \smsh{} Z)$,
so it suffices to prove that $|LM_.(X)| \to M(X)$ is a weak equivalence
for any $G$-CW complex $X$.  We will prove this by showing that it is
filtered by weak equivalences.

Suppose $X$ is a $G$-CW complex.  For any $n \in \mb N$, define
\[
M(X)^{(n)} = \left(\coprod_{|Y/G| \leq n, Y \in \Gamma^o_G} M(Y)
  \smsh{} F^G(Y,X)\right)/\sim.
\]
Here the equivalence relation is the same as that defining $M(X)$.

We now prove that $M(X)^{(n)}$ is a subspace of $M(X)$.
Suppose $Y \in \Gamma^o_G$, and $u \smsh{} v \in M(Y) \smsh{}
F^G(Y,X)$.  Call this element {\em minimal\/} if $v$ is an embedding
and $u$ is not in the image of $j_*$ for any proper inclusion $j_*$.

Given an element $u \smsh{} v$ as above, the map $v$ uniquely factors
through the surjection $p : Y \to {\rm Im}(v)$.  By
Remark~\ref{rmk:degen}, there is a unique minimal inclusion ${j : Z
\subset {\rm Im}(v)}$ such that $p_* u = j_* m$ for some $m$.  The
element $m \smsh{} j$ is minimal and equivalent to $u \smsh{} v$.
Note that if $|Y/G| \leq n$, this equivalence holds in $M(X)^{(n)}$.

If $f : Y' \to Y$, $u \in M(Y')$, and $v : Y \to X$, then there is a
natural factorization as follows:
\xym{
& Y' \ar[r]^f \ar@{->>}[d] & Y \ar@{>>}[d] \ar[dr]^v \\
Z\, \ar@{>->}[r] & {\rm Im}(f^* v)\, \ar@{>->}[r] & {\rm Im}(v)
\,\ar@{>->}[r]& X.
}
This shows that $f_* u \smsh{} v$ and $u \smsh{} f^* v$ are equivalent
to the same minimal element.

We find that two elements are equivalent under the equivalence
relation defining $M(X)$ or $M(X)^{(n)}$ if and only if they are
equivalent to a common minimal element.  As a result, an element in
$M(X)^{(n)}$ is not in the image of $M(X)^{(n-1)}$ if and only if it
is its own minimal factorization, i.e. it is of the form $u \smsh{} v$
where $v$ is an embedding and $u$ is not in the image of $j_*$ for $j$
any proper inclusion.  Two such minimal elements $u \smsh{} v$ and $u'
\smsh{} v'$ are equivalent if and only if there is an isomorphism $f$
such that $f_* u = u'$ and $(f^{-1})^* v = v'$.

In particular, for any $n > 0$, there is a natural pushout square for
constructing $M(X)^{(n)}$ from $M(X)^{(n-1)}$.  Define $\Sigma_n \wr
G$ to be the wreath product $G^n \rtimes \Sigma_n$, which is the
automorphism group in $\Gamma^o_G$ of $G_+ \smsh{} n_+$.  Let $\tilde
F(n_+, X)$ be the subset of $F(n_+,X) = F^G(G_+ \smsh{} n_+,X)$
consisting of those maps that are not embeddings, and let $\tilde
M(G_+ \smsh{} n_+)$ denote the union of the images of $M(Y)$ over
proper inclusions $Y \to G_+ \smsh{} n_+$.

There is a natural pushout diagram
\xym{
A \ar[r] \ar[d] & 
%\left(\tilde F(n_+,X)\smsh{\Sigma_n \wr G} M(G_+ \smsh{} n_+)\right)
%\cup \left(F^G(Z,X) \smsh{\Sigma_n \wr G} \tilde M(Z)\right)
M(X)^{(n-1)} \ar[d] \\
F(n_+,X) \smsh{\Sigma_n \wr G} M(G_+ \smsh{} n_+) \ar[r] &
M(X)^{(n)},
}where 
\[
A = \left(\tilde F(n_+,X)\smsh{\Sigma_n \wr G} M(G_+ \smsh{} n_+)\right)
\cup \left(F^G(Z,X) \smsh{\Sigma_n \wr G} \tilde M(G_+ \smsh{} n_+)\right).
\]
Similarly, there is a natural pushout diagram
\xym{
B \ar[r] \ar[d] & 
%\left(\tilde F(n_+,X)\smsh{\Sigma_n \wr G} M(G_+ \smsh{} n_+)\right)
%\cup \left(F^G(Z,X) \smsh{\Sigma_n \wr G} \tilde M(Z)\right)
|LM_.(X)|^{(n-1)} \ar[d] \\
F(n_+,X) \smsh{\Sigma_n \wr G} |LM_.(G_+ \smsh{} n_+)| \ar[r] &
|LM_.(X)|^{(n)},
}where $B$ is the corresponding union for $LM_.$.

Because $X$ is a free $G$-CW complex, the map $\tilde F(n_+,X) \to
F(n_+,X)$ is a cofibration of $(\Sigma_n \wr G)$-spaces.
Additionally, the map of $(\Sigma_n \wr G)$-spaces $|LM_.(G_+ \smsh{}
n_+)| \to M(G_+ \smsh{} n_+)$ is a weak equivalence, so the square
\xym{
\tilde F(n_+,X) \smsh{\Sigma_n \wr G} |LM_.(G_+ \smsh{} n_+)| \ar[r] \ar[d] &
\tilde F(n_+,X) \smsh{\Sigma_n \wr G} M(G_+ \smsh{} n_+) \ar[d] \\
F(n_+,X) \smsh{\Sigma_n \wr G} |LM_.(G_+ \smsh{} n_+)| \ar[r] &
F(n_+,X) \smsh{\Sigma_n \wr G} M(G_+ \smsh{} n_+)
}is homotopy cocartesian.  Similarly, the weak equivalence
$|\widetilde{LM}_.| \to \tilde M$ shows that the square
\xym{
\tilde F(n_+,X) \smsh{\Sigma_n \wr G} \left|\widetilde{LM}_.(G_+
  \smsh{} n_+)\right| \ar[r] \ar[d] &
\tilde F(n_+,X) \smsh{\Sigma_n \wr G} \tilde M(G_+ \smsh{} n_+) \ar[d] \\
F(n_+,X) \smsh{\Sigma_n \wr G} \left|\widetilde{LM}_.(G_+ \smsh{}
  n_+)\right| \ar[r] &
F(n_+,X) \smsh{\Sigma_n \wr G} \tilde M(G_+ \smsh{} n_+)
}is homotopy cocartesian.

These two previous homotopy cocartesian squares imply that the square
\xym{
B \ar[r] \ar[d] &
F(n_+,X) \smsh{\Sigma_n \wr G} |LM_.(G_+ \smsh{} n_+)| \ar[d] \\
A \ar[r] &
F(n_+,X) \smsh{\Sigma_n \wr G} M(G_+ \smsh{} n_+)
}is homotopy cocartesian.  (This square would be honestly cocartesian
if the previous squares were cocartesian.)  The horizontal maps in
this square are cofibrations.

Inductively assume that $|LM_.(X)|^{(n-1)} \to M(X)^{(n-1)}$ is a
weak equivalence.  Weak equivalences are preserved by pushouts along
cofibrations, so the map
\[
|LM_.(X)|^{(n-1)} \cup_{A} \left(F(n_+,X) \smsh{\Sigma_n \wr G} M(G_+
\smsh{} n_+)\right) \to M(X)^{(n)}
\]
is a weak equivalence.  However, we know that there is a weak
equivalence
\[
F(n_+,X) \smsh{\Sigma_n \wr G} M(G_+ \smsh{} n_+) \simeq
A \cup_B \left(F(n_+,X) \smsh{\Sigma_n \wr G} |LM_.(G_+ \smsh{} n_+)|\right),
\]
so we find that the map $|LM_.(X)|^{(n)} \to M(X)^{(n)}$ is a weak
equivalence, as desired.
\end{pf}

\begin{cor}
If $M$ is a semi-cofibrant $\Gamma_G$-space and a map $X \to Y$ of free
$G$-CW complexes is $k$-connected, so is the map $M(X) \to M(Y)$.
\end{cor}

\begin{pf}
It suffices to show that the map $|LM_.(X)| \to |LM_.(Y)|$ is
$k$-connected.  If $k = 0$, this is clear.

The simplicial space $LM_.(X)$ is good in the sense of
Segal~\cite{catcohom}, Appendix A, so it suffices to show that the map
of thick geometric realizations $||LM_.(X)|| \to ||LM_.(Y)||$ is
$k$-connected.  However, this map of simplicial spaces is
levelwise of the form

\xym{
\bigvee_{Z_0,\ldots,Z_p} F^G(Z_p,X) \smsh{} \Gm{Z_{p-1}}{Z_p} \smsh{}
\cdots \smsh{} \Gm{Z_0}{Z_1} \smsh{} M(Z_0) \ar[d]\\
\bigvee_{Z_0,\ldots,Z_p} F^G(Z_p,Y) \smsh{} \Gm{Z_{p-1}}{Z_p} \smsh{}
\cdots \smsh{} \Gm{Z_0}{Z_1} \smsh{} M(Z_0).
}

This map is $k$-connected because the map $F^G(Z_p,X) \to F^G(Z_p,Y)$
is.  The result follows because a levelwise $k$-connected map of
simplicial spaces $A_. \to B_.$ induces a $k$-connected map of thick
geometric realizations.  We include a proof as follows.

Filtering the thick geometric realization by skeleta, we find that
for any $n \geq 1$ there is a commutative square
\xym{
{\rm sk}_{n-1}(||A_.||) \ar[r] \ar[d] &
{\rm sk}_n(||A_.||) \ar[d] \\
{\rm sk}_{n-1}(||A_.||) \ar[r] &
{\rm sk}_n(||B_.||). \\
}
Assume inductively that
the leftmost vertical map is $k$-connected.

We now recall the statement of the Blakers-Massey excision theorem, as
in~\cite{goodwillie}, Section~2.  Suppose there is a commutative
square of spaces

\xym{
A \ar[r]^{f_1} \ar[d]_{f_2} & B \ar[d] \\ C \ar[r] & D.
}

If the map from the homotopy pushout to $D$ is $\ell$-connected, and
$f_i$ is $k_i$ connected for each $i$, then the map of of homotopy
fibers
\[
{\rm fib}(A \to C) \to {\rm fib}(B \to D)
\]
is ${\rm min}(k_1 + k_2 - 1,\ell - 1)$-connected.

Applying this shows that the map from $sk_{n-1} (||A_.||)$ to
$sk_n(||A_.||)$ is $(n-1)$-connected, and similarly for $||B_.||$.
Consider the commutative square of skeleta.  Because $k \geq 1$, the map
from the homotopy pushout to $sk_n(||B_.||)$ is $1$-connected by
the Seifert-Van Kampen theorem.  By the relative Hurewicz theorem, the 
connectivity of this map is the same as the connectivity of its
homotopy cofiber.

The total homotopy cofiber of the square of skeleta is the
cofiber of the map $S^n \smsh{} LM_n(X)_+ \to S^n \smsh{} LM_n(Y)_+$,
which is $(n+k)$-connected by assumption.  The map of $(n-1)$-skeleta
is $k$-connected, and the map from the $(n-1)$-skeleton to the
$n$-skeleton is $(n-1)$-connected.  Therefore, the Blakers-Massey
theorem shows that the map of homotopy fibers is $(n + k -
2)$-connected.  As $n \geq 1$, it is in particular
$(k-1)$-connected.  The homotopy fiber of the map of $(n-1)$-skeleta
is $(k-1)$-connected, so the homotopy fiber of the map of $n$-skeleta
must be $(k-1)$-connected as well, as desired.
\end{pf}

For any $\Gamma_G$-space $M$, we have an associated (na\"ive
pre-)spectrum $\{M(G_+ \smsh{} S^n)\}$, which is the spectrum
associated to the underlying $\Gamma$-space $M(G_+\smsh{} -)$.  A map
of $\Gamma_G$-spaces $M \to M'$ is called a stable equivalence if the
associated map of spectra is a weak equivalence.

\begin{prop}
\label{prop:smash}
For any semi-cofibrant $\Gamma_G$-space $M$ and free based $G$-CW complex
$X$, the assembly map $X \smsh{G} M(G_+ \smsh{} -) \to X \otimes_G M$ is a
stable equivalence.
\end{prop}

\begin{pf}
It suffices to show that $X \smsh{G} M(G_+ \smsh{} S^n) \to
M(X \smsh{} S^n)$ is highly connected for large $n$.
Using the levelwise weak equivalence $|X \otimes_G LM_.| \to X
\otimes_G M$, it suffices to show that this statement is true for
$\Gamma$-spaces of the form $\Gm{Y}{-}$ for $Y \in \Gamma^o_G$.

In this case, we have the following diagram:

\xym{
X \smsh{G} \bigvee_Y (G_+ \smsh{} S^n) \ar[r] \ar[d] &
\bigvee_Y X \smsh{G} (G_+ \smsh{} S^n) \ar[d] \\
X \smsh{G} \prod_Y (G_+ \smsh{} S^n) \ar[r] \ar@{=}[d]&
\prod_Y (X \smsh{} S^n) \ar@{=}[d]\\
X \smsh{G} \Gm{Y}{G_+ \smsh{} S^n} \ar[r]&
\Gm{Y}{X \smsh{} S^n}.
}

The top vertical arrows are isomorphisms on homotopy groups up to roughly
dimension $2n$, as $G_+ \smsh{} S^n$ is $(n-1)$-connected.  The
uppermost horizontal arrow is an isomorphism.  Therefore, the bottom
map is an equivalence on homotopy groups up to roughly dimension $2n$,
as desired.
\end{pf}

\section{Connective equivariant $K$-homology}
\label{sec:ctkh}

In this section, we construct for each $n$ a $\Gamma_{\PU(n)}$-space
whose underlying spectrum is homotopy equivalent to $ku$, the
connective $K$-theory spectrum.

Fix an integer $n$.  For any $d \in \mb N$, we have the Stiefel
manifold $V(nd)$ of isometric embeddings of $\mb C^n \otimes \mb C^d$
into $\mb C^\infty$, where these vector spaces have the standard inner
products.  The tensor product map $\U(n) \otimes \U(d) \to \U(nd)$
gives the space $V(nd)$ a free right action of $I \otimes \U(d)$ by
precomposition, where $I$ is the identity element of $\U(n)$.  Denote
the quotient space by $H(d)$, and write $H = \coprod_d H(d)$.  (Note
$H \simeq \coprod_d \BU(d)$.)  The space $H$ has a partially defined
direct sum operation: if $\{V_i\}$ is a finite set of elements of $H$
such that $V_i \perp V_j$ for $i \neq j$, there is a sum element $\oplus
V_i$ in $H$.

There is also an action of $\U(n) \otimes I$ on $V(nd)$ that commutes
with the action of $I \otimes \U(d)$, and hence passes to an action on
the quotient $H(d)$.  Because $\lambda I \otimes I = I \otimes \lambda
I$, the scalars in $\U(n)$ act trivially on $H(d)$, so the action
factors through $\PU(n)$.  We therefore get a right action of $\PU(n)$
on $H$.  The direct sum operation is $\PU(n)$-equivariant.

For any $Z \in \Gamma^o_{\PU(n)}$, define
\[
ku^{\PU(n)}(Z) = \left\{ f \in F^{\PU(n)}(Z,H)\ \Big|\
  f(z) \perp f(z') \hbox{ if } [z] \neq [z'] \right\}.
\]
A point of $ku^{\PU(n)}(Z)$ consists of a subspace $\mb C^\infty$,
isomorphic to $\mb C^{{\rm dim} f(z)}$, associated to each 
non-basepoint $[z]$ of $Z/\PU(n)$, such that the vector spaces
associated to $[z]$ and $[z']$ are orthogonal if $[z] \neq [z']$.

Given a map $\alpha \in \Gamma^o_{\PU(n)}(Z,Z')$ and $f \in
ku^{\PU(n)}(Z)$, we get an element $ku^{\PU(n)}(\alpha)(f) \in
ku^{\PU(n)}(Z')$ as follows:
\[
ku^{\PU(n)}(\alpha)(f)(z') = \bigoplus_{\alpha(z) = z'} f(z).
\]
This direct sum is well-defined: if the preimage of $z$ is the family
$\{z_i\}$, then the $z_i$ all lie in distinct orbits because the
action of $\PU(n)$ is free away from the basepoint.  Therefore, the
subspaces associated to the $z_i$ are orthogonal.  The map
$ku^{\PU(n)}(\alpha)(f)$ is also clearly $\PU(n)$-equivariant, and
takes distinct orbits to orthogonal elements of $H$.  

The underlying $\Gamma$-space is given as follows.  For $Z \in
\Gamma^o$,
\[
ku^{\PU(n)}(\PU(n)_+ \smsh{} Z) = \{f \in F(Z,H)\ | f(z) \perp f(z')
\hbox{ if } z \neq z'\}.
\]
The spectrum attached to the underlying $\Gamma$-space of
$ku^{\PU(n)}$ is weakly equivalent to the connective $K$-theory
spectrum $ku$---see~\cite{segal}.

The $\Gamma_{\PU(n)}$-space $ku^{\PU(n)}$ is semi-cofibrant: The image in 
$ku^{\PU(n)}(Z)$ of the spaces $ku^{\PU(n)}(Y)$ for $Y \subsetneq Z$
consists of those maps $Z \to H$ that map some nontrivial subset of
$Z$ to $H(0)$.  This is a union of components of $ku^{\PU(n)}(Z)$.

We also define a second $\Gamma_{\PU(n)}$-space $ku/\beta$ as follows.  
For any $z \in \Gamma^o_{PU(n)}$,
\[
ku/\beta(Z) = \tilde{\mb N}[Z/{\PU(n)}],
\]
where $\tilde{\mb N}$ is the reduced free abelian monoid functor.
More explicitly, $ku/\beta(Z)$ is the quotient of the free abelian
monoid on $Z/{\PU(n)}$ by the submonoid $\mb N [*]$.  For $\alpha \in
\Gamma^o_{\PU(n)}(Z,Z')$,
\[
ku/\beta(\alpha)\left(\sum n_z[z]\right) = \sum n_z[\alpha(z)].
\]
(The reason for the notation is that the underlying spectrum is the
cofiber of the Bott map.)  The $\Gamma_{\PU(n)}$-space $ku/\beta$ is
semi-cofibrant for the same reason as $ku^{\PU(n)}$.

For $X$ a free right $\PU(n)$-space, $X \otimes_{\PU(n)} ku/\beta$ is
the infinite symmetric product ${\rm Sym}^\infty(X/\PU(n))$.

There is a natural map $\epsilon: ku^{\PU(n)} \to ku/\beta$ of
$\Gamma_{\PU(n)}$-spaces: if $f \in ku^{\PU(n)}(Z)$, define $\epsilon(f) =
\sum_{[z]} \left(\frac{\dim f(z)}{n}\right) [z]$.  
The map $\epsilon$ represents the augmentation $ku \to \eilm{\mb Z}$
on the underlying spectra; it is the first stage of the Postnikov
tower for $ku$.

The $\Gamma_{\PU(n)}$-space $ku^{\PU(n)}$ determines a homology theory for
$\PU(n)$-spaces.  Specifically, we can define
\[
ku^{\PU(n)}_*(X) = \pi_*\left(X \otimes_{\PU(n)} ku^{\PU(n)}\right).
\]
Here $\pi_*$ denotes the stable homotopy groups of the spectrum.  In
fact, because the underlying spectrum of $ku^{\PU(n)}$ is {\em special},
we can compute 
\[
ku^{\PU(n)}_*(X) = \pi_*\left(ku^{\PU(n)}(X)\right)
\]
for $X$ connected.  (See~\cite{catcohom}, 1.4.)

\section{Unitary deformation $K$-theory}
\label{sec:defk}

In this section we will have a fixed finitely generated {\em discrete}
group ${\rm G}$.  Carlsson, in~\cite{carlsson:kth}, defined a notion
of the ``deformation $K$-theory'' of ${\rm G}$ as a contravariant
functor from groups to spectra, and in the introduction of this
article an analogous notion of ``unitary deformation $K$-theory''
${\cal K}{\rm G}$ was sketched.  The following are weakly equivalent
definitions of the corresponding notion of ${\cal K}{\rm G}$:
\begin{itemize}
\item The spectrum associated to the the $E_\infty$-$H$-space
\[
\coprod_n \EU(n) \pull{\U(n)} \Hom({\rm G}, \U(n)).
\]
\item The $K$-theory of a category of unitary representations of
  ${\rm G}$.  This category is an internal category in \cat{Top}:
  i.e., the objects and morphism sets are both given topologies.
\item The $K$-theory of the singular complex of the category above.
  (This is essentially the definition given in~\cite{carlsson:kth}.)
\end{itemize}

We will now describe another model for the unitary deformation
$K$-theory of ${\rm G}$, equivalent to the first definition given
above.  The construction is based on the construction of connective
topological $K$-homology of Segal in~\cite{segal}.  (Also
see~\cite{walker}.)

\begin{defin}
Let ${\mathcal U} = \mb C^\infty$ be the infinite inner product space
having orthonormal basis $\{e_i\}_{i=0}^{\infty}$, with action of the
group $\U = \colim \U(n)$.
\end{defin}

\begin{defin}
A ${\rm G}$-plane $V$ of dimension $k$ is a pair
$(V,\rho)$, where $V$ is a $k$-dimensional plane in ${\mathcal U}$ and
$\rho : {\rm G} \to \U(V)$ is an action of ${\rm G}$ on $V$.
\end{defin}

We now describe a (non-equivariant) $\Gamma$-space ${\cal K}{\rm G}$.  Define
\[
{\cal K}{\rm G}(X) = \left\{(V_x,\rho_x)_{x \in X}\ \Big |\ V_x\hbox{ a
  }{\rm G}\hbox{-plane, }V_x \perp V_{x'}\hbox{ if } x \neq x', V_* =
  0 \right\}.
\]
This is a special $\Gamma$-space.  The underlying $H$-space is
\[
{\cal K}{\rm G}(1_+) \cong \coprod_n V(n) \times_{\U(n)} \Hom({\rm G},
\U(n)),
\]
where $V(n)$ is the Stiefel manifold of $n$-frames in ${\cal
U}$.  We will now describe the simplicial space $X_. = {\cal
K}{\rm G}(S^1)$.  Because ${\cal K}{\rm G}$ is special, $\Omega |X_.| \simeq
\Omega^\infty {\cal K}{\rm G}$.

For $p > 0$, $X_p$ is the space
\[
\left\{(V_i,\rho_i)_{i=1}^{p}\ \Big |\ (V_i,
  \rho_i) \hbox{ a ${\rm G}$-plane, } V_i \perp V_j \hbox{ if } i \neq
  j \right\}.
\]
($X_0$ is a point.)  Face maps are given by taking sums
of orthogonal ${\rm G}$-planes or removing the first or last
${\rm G}$-plane.  Degeneracy maps are given by insertion of
$0$-dimensional ${\rm G}$-planes.

The geometric realization of this simplicial space can be explicitly
identified.  Let $Y$ be the space of pairs $(A,\rho)$, where $A \in
\U$ and $\rho$ is a homomorphism ${\rm G} \to \U$ such that $\rho(g)A
= A \rho(g)$ for all $g \in {\rm G}$.  Call two such elements $(A,
\rho)$ and $(A', \rho')$ equivalent if $A = A'$ and $\rho, \rho'$
agree on all eigenspaces of $A$ corresponding to eigenvalues $\lambda
\neq 1$.  Write the standard $p$-simplex $\Delta^p$ as the set of all
$0 \leq t_1 \leq \ldots \leq t_p \leq 1$. Then there is a
homeomorphism $|X_.| \to (Y\ /\sim)$ given by sending a point $((V_i,
\rho_i)_{i=1}^{p}, 0 \leq t_1 \leq \ldots \leq t_p \leq 1)$ of $X_p
\times \Delta^p$ to the pair $(A,\rho)$, where $A$ acts on $V_i$ with
eigenvalue $e^{2\pi i t_i}$ and by $1$ on the orthogonal complement of
$\Sigma V_i$, while $\rho$ acts on $V_i$ by $\rho_i$ and acts by $1$
on the orthogonal complement of $\Sigma V_i$.  (Here $\Sigma V_i$ is
the span of the set of orthogonal subspaces $V_i$.)  This map is a
homeomorphism by the spectral theorem.  (The essential details of this
argument are from~\cite{harris} and~\cite{mitchell}.)

We will refer to this space $|X.| \cong (Y\ /\sim)$ as $E$.  It is
space $1$ of the $\Omega$-spectrum associated to ${\mathcal K}{\rm
G}$, in the sense that $\Omega^\infty{\cal K}{\rm G} \simeq \Omega E$.

This method is applicable to various other categories of
representations of ${\rm G}$ that we will now examine in detail.

For any $n \geq 0$, there is a sub-$\Gamma$-space ${\cal K}{\rm G}_n$ of
${\cal K}{\rm G}$.  The space ${\cal K}{\rm G}_n(X)$ consists of those elements
$\{(V_x,\rho_x)\}_{x \in X}$ of ${\cal K}{\rm G}(X)$ such that $\rho_x$ breaks
up into a direct sum of irreducible representations of dimension less
than or equal to $n$.  Each ${\cal K}{\rm G}_n$ is a special $\Gamma$-space.

\begin{prop}
\label{prop:hocolim}
The map $\hocolim {\cal K}{\rm G}_n \to {\cal K}{\rm G}$ is a weak equivalence.
\end{prop}

\begin{pf}
Clearly ${\cal K}{\rm G}$ is the union of the sub-$\Gamma$-spaces ${\cal
K}{\rm G}_n$.  For any based set $X$, any element $\{(V_x,\rho_x)\}$ of
${\cal K}{\rm G}(X)$ has a well-defined total dimension $\sum \dim V_x$,
and this dimension is locally constant.  Therefore, ${\cal K}{\rm G}(X)$
breaks up as a disjoint union according to total dimension.  The
component consisting of elements of total dimension $n$ is completely
contained in the subspace ${\cal K}{\rm G}_N(X)$ for all $N \geq n$.  In
particular, if $x$ is any point of ${\cal K}{\rm G}(X)$ whose total
dimension is $n$, $\pi_*({\cal K}{\rm G}(X),x) \cong \pi_*({\cal
K}{\rm G}_N(X),x)$ for all $N \geq n$.
\end{pf}

\begin{rmk}
As a result, for $X$ a simplicial set finite in each dimension (such as
a sphere), $\hocolim {\cal K}{\rm G}_n(X)$ can be formed levelwise, and is
levelwise weakly equivalent to ${\cal K}{\rm G}(X)$.  In particular, the
natural map $\hocolim {\cal K}{\rm G}_n(S^k) \to {\cal K}{\rm G}(S^k)$ is a weak
equivalence for all $k$, so the associated spectrum of ${\cal K}{\rm G}$ is
weakly equivalent to the associated spectrum of $\hocolim {\cal
K}{\rm G}_n$.
\end{rmk}

We have infinite loop spaces $E_n = |{\cal K}{\rm G}_n(S^1)|$.  For any $n
\in \mb N$, $E_n$ is the subspace of $E$ consisting of pairs $(A,
\rho)$ such that $\rho$ is a direct sum of irreducible representations
of ${\rm G}$ of dimension less than or equal to $n$.

This gives a sequence of inclusions
\[
* = E_0 \subset E_1 \subset E_2 \subset \ldots
\]
of infinite loop spaces.  Each of these inclusions is part of a
quasifibration sequence $E_{n-1} \to E_n \to B_n$ where the base
spaces will be explicitly identified.  This gives rise to the
following ``unrolled exact couple'' of infinite loop spaces:

\xym{
{*} \ar[rr] && E_1 \ar[dl] \ar[rr] && E_2 \ar[dl] \ar[rr] &&
E_3\ldots \ar[dl]\\
& B_1 \ar@{>}|O[ul] && B_2 \ar@{>}|O[ul] && \ar@{>}|O[ul]
B_3.
}

Additionally, the inclusions of infinite loop spaces $E_n$ are induced
by maps of $E_\infty$-$ku$-modules, so the above is induced by an exact
couple of $ku$-module spectra.

The intuition for the description of $B_n$ is that the category of
${\rm G}$-representations whose irreducible summands have dimension
less than $n$ forms a Serre subcategory of the category of ${\rm
G}$-representations whose irreducibles have dimension less than or
equal to $n$, and the quotient category should be the category of sums
of irreducible representations of dimension exactly $n$.  The topology
on the categories involved complicates the question of when a
localization sequence of spectra exists in this situation, as the
most obvious attempts to generalize of Quillen's Theorem B would not
be applicable.  We will construct the localization sequence
explicitly.

There is a quotient $\Gamma$-space $F_n$ of ${\cal K}{\rm G}_n$ by the
equivalence relation $(V_x, \rho_x)_{x \in X} \sim (V'_x, \rho'_x)_{x
\in X}$ if for all $x \in X$:
\begin{itemize}
\item The subspace $W_x$ of $V_x$ generated by irreducible
  subrepresentations of $\rho$ of dimension $n$ coincides with the
  corresponding subspace for $\rho'$, and
\item $\rho$ and $\rho'$ agree on $W_x$.
\end{itemize}
Again, $F_n$ is a special $\Gamma$-space.  

Define $B_n$ to be the space $|F_n(S^1)|$.  $B_n$ is the quotient of
$E_n$ by the following equivalence relation.  We say $(A,\rho) \sim
(A', \rho')$ if:
\begin{itemize}
\item The subspace $W$ of $\mathcal{U}$ generated by irreducible
  subrepresentations of $\rho$ of dimension $n$ is the same as the
  subspace of ${\mathcal U}$ generated by irreducible
  subrepresentations of $\rho'$ of dimension $n$,
\item $\rho$ and $\rho'$ have the same action on $W$, and
\item $A$ and $A'$ have the same action on $W$.
\end{itemize}

Note that each equivalence class contains a unique pair $(A,\rho)$
such that $\rho$ acts trivially on the eigenspace of $A$ for $1$, and
on the complementary subspace $\rho$ is a direct sum of irreducible
$n$-dimensional representations.

\section{Proof of the existence of the localization sequence}
\label{sec:local}

The proof that $p_n: E_n \to B_n$ is a quasifibration (and hence
induces a long exact sequence on homotopy groups) proceeds inductively
using the following result.

\begin{thm}[Hardie~\cite{hardie}]
Suppose that we have a diagram

\xym{
Q \ar[d]^\lambda & f^*(E) \ar[l]_h \ar[r] \ar[d]^s & E \ar[d]^p \\
Q' & A \ar[l]^g \ar[r]_f & B.
}

Here $f$ is a cofibration, $p$ is a fibration, $f^*(E)$ is the
pullback fibration, and $\lambda$ is a quasifibration.  If $h :
s^{-1}(a) \to \lambda^{-1}(ga)$ is a weak equivalence for all $a \in
A$, then the induced map of pushouts $Q \coprod_{f^*(E)} E \to Q' \coprod_A B
$ is a quasifibration.
\end{thm}

\begin{prop}
The map $p_n: E_n \to B_n$ is a quasifibration with fiber
$E_{n-1}$.
\end{prop}

\begin{rmk}
This is what we might expect, as the map $E_n \to B_n$ is precisely
the map that forgets the irreducible subrepresentations of dimension
less than $n$.  The fact that $E_{n-1}$ is the honest fiber over any
point is clear, but we need to show that $E_{n-1}$ is also the homotopy
fiber.
\end{rmk}

\begin{pf}
We will proceed by making use of a rank filtration.  These rank
filtrations were introduced in~\cite{mitchell} and~\cite{rognes}.
In particular, Mitchell explicitly describes this rank filtration
for the connective $K$-theory spectrum.

For any $j$, let $B_{n,j}$ be the subspace of $B_n$ generated by those
pairs $(A,\rho)$ such that $\rho$ contains at most a sum of $j$
irreducible representations.  There is a sequence of inclusions
$B_{n,j-1} \subset B_{n,j}$.  Write $E_{n,j}$ for the subset of $E_n$
lying over $B_{n,j}$.

The map $E_{n,0} \to B_{n,0}$ is a quasifibration, because $B_{n,0}$
is a point. Now suppose inductively that $E_{n,j-1} \to B_{n,j-1}$ is
a quasifibration.

Let $Y_j$ be the space of triples $(A,\rho,W)$, where $W$ is an
$nj$-dimensional subspace of ${\mathcal U}$, $A$ is an element of
$\U(W)$, and $\rho$ is a representation of ${\rm G}$ on $W$ commuting
with $A$ and containing irreducible summands of dimension $n$ or less.
Let $X_j$ be the subset of $Y_j$ of triples $(A,\rho,W)$ such that
$(A, \rho)$ represents a pair in $B_{n,j-1}$; in other words, $\rho$
contains less than $j$ distinct $n$-dimensional irreducible summands
on the orthogonal complement of the eigenspace for $1$ of $A$.

Next, we define a space $Y'_j$ of triples $(A, \rho, W)$, where
$(A,\rho) \in E_{i,j}$ and $W$ is an $A$- and $\rho$-invariant
$nj$-dimensional subspace of ${\mathcal U}$ containing all
the $n$-dimensional irreducible summands of $\rho$.  There is a map $Y'_j
\to Y_j$ given by forgetting the actions of $A$ and $\rho$ off $W$.
Let $X'_j$ be the fiber product of $X_j$ and $Y'_j$ over $Y_j$; it
consists of triples $(A, \rho, W)$ where $\rho$ contains less than $j$
distinct $n$-dimensional summands.

There is a map $X_j \to B_{j-1}$ given by sending $(A, \rho, W)$ to
$(A, \rho)$, and a similar map $X'_j \to E_{j-1}$.  These maps all
assemble into the diagram below.

\xym{
E_{n,j-1} \ar[d]^{p_n} & X'_j \ar[l] \ar[r] \ar[d] & Y'_j \ar[d]^p \\
B_{n,j-1} & X_j \ar[l] \ar[r] & Y_j
}

There is an evident map from the pushout of the bottom row to
$B_{n,j}$, and similarly a map from the pushout of the top row to
$E_{n,j}$.

The map $X_j \to B_{n,j-1}$ is a quotient map; two points become
identified by forgetting the ``framing'' subspace $W$, the
non-$n$-dimensional summands of $\rho$, and the summands of $\rho$ on
the eigenspace for $1$ of $A$.  For points of $Y_j$ not in $X_j$, the
framing subspace $W$ is determined by the image $(A,\rho)$ in $B_j$
because $\rho$ must have $j$ distinct $n$-dimensional irreducible
summands covering all of $W$, and $A$ can have no eigenspace for the
eigenvalue $1$.  Therefore, the map from $Y_j$ to the pushout of the
bottom row is precisely the quotient map gotten by forgetting the
framing $W$ and any non-$n$-dimensional summands or summands lying on
the eigenspace for $1$ of $A$.  This identifies the pushout
with $B_{n,j}$.  In exactly the same way, the pushout of the top row
is $E_{n,j}$.  The induced map of pushouts is the projection map
$E_{n,j} \to B_{n,j}$.

The map $X_j \to Y_j$ is a cofibration because it is the colimit of
geometric realizations of a closed inclusion of real points of
algebraic varieties. (The subspace $W$ is allowed to vary over the
infinite Grassmannian.  If we restrict its image to any finite
subspace we get an inclusion of real algebraic varieties.)

The right-hand square is a pullback by construction, and the map $p_n$
is assumed to be a quasifibration.

The map $p$ is a fiber bundle with fiber $E_{n-1}$: An equivalence
class of points $(A,\rho,W) \in Y'_j$ consists of a choice of
$nj$-dimensional subspace $W$ of ${\mathcal U}$, a choice of element
in $(\bar A, \bar \rho, W)$ in $Y_j$ to determine the action of $A$
and $\rho$ on $W$, and a choice of $(A', \rho')$ acting on the
orthogonal complement of $W$ such that $\rho'$ is made up of summands
of dimension less than $n$.  In other words, there is a pullback
square:

\xym{
Y'_j \ar[r] \ar[d] & V\ar[d] \\
Y_j \ar[r] & {\rm Gr}(nj).
}

Here ${\rm Gr}(nj)$ is the Grassmannian of $nj$-dimensional planes in
${\mathcal U}$, and $V$ is the bundle over the Grassmannian consisting
of $nj$-dimensional planes in $\mathcal{U}$ and elements of
$E_{n-1}$ acting on their orthogonal complements.

Given any point $(A,\rho,W)$ of $X_j$, the fiber in $X_j'$ is
$E_{n-1}$ acting on the orthogonal complement of $W$.  Suppose that
$(A,\rho)$ in $B_{n,j-1}$ is in canonical form: $\rho$ acts by
a sum of irreducible dimension $n$ representations on some subspace
$W' \subset W$ and trivial representations on the orthogonal
complement, and $A$ has eigenvalue $1$ on the orthogonal complement of
$W'$.  Then the fiber over $(A,\rho)$ in $E_{n,j-1}$ consists of all
possible actions of $E_{n-1}$ on the orthogonal complement of $W'$.
The map from the fiber over $(A,\rho,W)$ to the fiber over $(A,\rho)$
is the inclusion of $E_{n-1}$ acting on $W^{\perp}$ to $E_{n-1}$
acting on $(W')^{\perp}$.  This inclusion is a homotopy equivalence.

Therefore, $E_{n,j} \to B_{n,j}$ is a quasifibration with fiber
$E_{n-1}$.  Taking colimits in $j$, $E_n \to B_n$ is a quasifibration
with fiber $E_{n-1}$.
\end{pf}

\begin{cor}
\label{cor:fibseq}
The maps ${\cal K}{\rm G}_{n-1} \to {\cal K}{\rm G}_n \to F_n$ realize to a
fibration sequence in the homotopy category of spectra.
\end{cor}

\begin{pf}
This follows because the composite map is null, the map from ${\cal
K}{\rm G}_{n-1}(S^1)$ to the homotopy fiber of ${\cal K}{\rm G}_n(S^1) \to
F_n(S^1)$ is a weak equivalence, and all three of these
$\Gamma$-spaces are special.
\end{pf}

\section{Identification of the $\Gamma$-space $F_n$}
\label{sec:identify}

Using the results of section~\ref{sec:gamma}, we will now identify the
$\Gamma$-spaces $F_n$ as equivariant smash products.

Let ${\rm Sum}({\rm G},n)$ be the subspace of $\Hom({\rm G}, \U(n))$
of reducible ${\rm G}$-representations of dimension $n$.  Define $R_n
= \Hom({\rm G}, \U(n)) / {\rm Sum}({\rm G},n)$.  There is a free
right action of $\PU(n)$ on $R_n$ by conjugation.

\begin{rmk}
The action of $\PU(n)$ on $R_n$ is free because the only
endomorphisms of an irreducible complex representation are scalar
multiplications.  This fails for orthogonal representations.
\end{rmk}

According to a result of Park and Suh~(\cite{park02}, Theorem 3.7),
the space $\Hom({\rm G},U(n))$, which is the set of real points of an
algebraic variety, admits the structure of a $\U(n)$-CW complex.  All
isotropy groups contain the diagonal subgroup, so this structure is
actually the structure of a $\PU(n)$-CW complex.  The subspace ${\rm
Sum}({\rm G},n)$ consists of those elements of $\Hom({\rm
G},\U(n))$ that are not acted on freely by $\PU(n)$, and so it must
be a CW-subcomplex.  Therefore, $R_n$ has an induced CW-structure.

\begin{prop}
\label{prop:baseid}
There is an isomorphism of $\Gamma$-spaces
\[
R_n \otimes_{\PU(n)} ku^{\PU(n)} \to F_n.
\]
\end{prop}

\vskip 0.7pc
\begin{pf}
By the universal property of the coend
\[
R_n \otimes_{\PU(n)} ku^{\PU(n)}(Z) = \int^Y ku^{\PU(n)}(Y) \smsh{}
F^G(Y,R_n\smsh{}Z),
\]
we can construct the map by exhibiting maps
\[
ku^{\PU(n)}(Y) \smsh{} F^G(Y,R_n\smsh{}Z) \to F_n(Z),
\]
natural in $Z$, that satisfy appropriate compatibility relations in $Y$.

Recall that a point of $ku^{\PU(n)}(Y)$ consists of an equivariant map
$f: Y \to H = \coprod V(nd)/I \otimes \U(d)$ such that $f(y) \perp
f(y')$ if $y \neq y'$.

Suppose $f \smsh{} h \in ku^{\PU(n)}(Y) \smsh{} F^G(Y,R_n \smsh{} Z)$.
For every $y \in Y$ the element $h(y) = r(y) \smsh{} z(y)$ determines
an irreducible action $r(y)$ of ${\rm G}$ on $\mb C^n$.  The element
$f(y) \in V(nd)/I \otimes \U(d)$ is the image of some element
$\widetilde {f(y)} \in V(nd)$, which determines an isometric embedding
$\mb C^n \otimes \mb C^d \to {\cal U}$.  Combining these two gives an
action of ${\rm G}$ on an $nd$-plane of ${\cal U}$, together with a
marking $z(y)$ of the plane by an element of $Z$.  The action of $I
\otimes \U(d)$ commutes with the ${\rm G}$-action on $\mb C^n \otimes \mb
C^d$, so the choice of lift $\widetilde{f(y)}$ does not change the
resulting ${\rm G}$-plane.  For $g \in {\rm G}$, $r(yg) = r(y) \cdot
g = g^{-1} r(y) g$, and $f(yg) = (g^{-1} \otimes I) f(y) (g \otimes
I)$, so the resulting plane only depends on the orbit $yG$.  The
resulting ${\rm G}$-plane breaks up into irreducible summands of
dimension precisely $n$.  Assembling these ${\rm G}$-planes over the
distinct orbits gives a collection of orthogonal hyperplanes with
${\rm G}$-actions, marked by points of $Z$, that break up into a
direct sum of $n$-dimensional irreducible representations.  As $r(y)$
approaches the basepoint of $R_n$, the representation becomes
reducible, so the map determines a well-defined element of $F_n(Z)$.
The compatibility of this map with maps in $Y$ is due to the fact that
it preserves direct sums.

This map is bijective; associated to any point of $F_n(Z)$ there is a
unique equivalence class of points that map to it.

We now construct the inverse map $F_n \to R_n \otimes_{\PU(n)} ku^{\PU(n)}$.
As $F_n$ is a quotient of ${\cal K}{\rm G}_n$, it suffices to construct a
map 
\[
{\cal K}{\rm G}_n(Z) \to R_n \otimes_{\PU(n)} ku^{\PU(n)}(Z),
\]
natural in $Z$, that respects the equivalence relation.

For $d \in \mb N$, consider the space $X_d$ of pairs $(\rho,
\{e_i\})$ consisting of an action $\rho$ of ${\rm G}$ on $\mb C^d$,
together with a set of nonzero mutually orthogonal ${\rm
G}$-equivariant projection operators $e_i$ whose image each have
dimension $n$ and that contain all of the $n$-dimension summands of
$\mb C^d$.  This space is compact Hausdorff.  The forgetful map
to $\Hom({\rm G},\U(d))$ identifies the image subspace $X'_d$,
consisting of elements $\rho$ that are a direct sum of representations
of dimension $n$ or less, with the quotient of $X_d$ by the
equivalence relation gotten by forgetting the $e_i$.

Recall that $V(d)$ is the space of embeddings of $\mb C^d$ in ${\cal
U}$. We now define a natural map
\begin{eqnarray*}
\phi: V(d) &\times &X_d \to \\
&&\bigvee_m ku^{\PU(n)}(\PU(n)_+ \smsh{} m_+) \smsh{}
F^{\PU(n)}(\PU(n)_+\smsh{} m_+, R_n).
\end{eqnarray*}
Given an isometric embedding $f$, an action $\rho$ of ${\rm G}$ on
$\mb C^d$, and a collection of $m$ ${\rm G}$-planes of $\mb C^d$, we
use the isometric embedding to obtain a collection $\{f(\hbox{Im
}e_i)\}$ of $n$-planes in ${\cal U}$, marked by points of $R_n$.
Composing with the quotient map, we get a map
\[
\phi': V(d) \times X_d \to R_n \otimes_{\PU(n)} ku^{\PU(n)}(1_+).
\]
The image of a point depends only on the $n$-dimensional summands of
the representation $\rho$.

The map $\phi'$ is invariant under the choice of direct sum
decomposition $\{e_i\}$, so there is an induced map $\phi''$ from the
quotient space $V(d) \times X'_d$.  (This identification of the
quotient space requires that we use the product in the category of
compactly generated spaces.)  The map $\phi''$ is invariant under the
diagonal action of $\U(d)$ on $V(d) \times X'_d$, so there is an
induced map $V(d) \times_{\U(d)} X'_d \to R_n \otimes_{\PU(n)}
ku^{\PU(n)}(1_+)$. Assembling these together over $d$ gives a map
${\cal K}{\rm G}_n(1_+) \to R_n \otimes_{\PU(n)} ku^{\PU(n)}(1_+)$.

The map
\[
R_n \otimes_{\PU(n)} ku^{\PU(n)}(Z) \to \prod_Z R_n \otimes_{\PU(n)}
ku^{\PU(n)}(1_+)
\]
is an inclusion of the subspace of mutually orthogonal elements.  The
composite 
\[
{\cal K}{\rm G}_n(Z) \to \prod_Z {\cal K}{\rm G}_n(1_+) \to \prod_Z R_n \otimes_{\PU(n)}
ku^{\PU(n)}(1_+)
\]
maps into this subspace.  Therefore, this gives a lift to a map of
$\Gamma$-spaces ${\cal K}{\rm G}_n \to R_n \otimes_{\PU(n)} ku^{\PU(n)}$.

This map respects the equivalence relation defining $F_n$, as the
image of a point of ${\cal K}{\rm G}_n$ only depended on its $n$-dimensional
summands.
\end{pf}

\begin{cor}
\label{cor:weaksmash}
There is a stable equivalence of $\Gamma$-spaces
\[
R_n \smsh{\PU(n)} ku^{\PU(n)} \to F_n.
\]
\hbox{}
\end{cor}

\begin{pf}
This follows from Proposition~\ref{prop:smash}.
\end{pf}

\section{$E_\infty$-algebra and module structures}
\label{sec:einfty}

In this section we will make explicit the following.
The tensor product of representations leads to the following
multiplicative structures:
\begin{itemize}
\item $ku$ is an $E_\infty$-ring spectrum,
\item ${\cal K}{\rm G}$ is an $E_\infty$-algebra over $ku$,
\item the sequence of maps ${\cal K}{\rm G}_1 \to {\cal K}{\rm G}_2 \to
  \cdots \to {\cal K}{\rm G}$ is a sequence of $E_\infty$-$ku$-module maps,
\item there are compatible $E_\infty$-$ku$-linear pairings ${\cal
    K}{\rm G}_n \smsh{} {\cal K}{\rm G}_m \to {\cal K}{\rm G}_{nm}$
  for all $n,m$, and
\item the assembly map $R_n \smsh{\PU(n)} ku^{\PU(n)} \to {\cal K}{\rm G}_n$
  is a map of $E_\infty$-$ku$-modules.
\end{itemize}
All of the above structures are natural in ${\rm G}$.

To begin, we will first recall the definition of a {\em
multicategory\/}.  A multicategory is an ``operad with several
objects'', as follows.  See~\cite{rma}.

\begin{defin}
A multicategory $\multc M$ consists of the following data:
\begin{itemize}
\item a class of objects $\Ob(\multc M),$
\item a set $\multc M_k(a_1,\ldots,a_k;b)$ for each $a_1,\ldots,a_k,b \in
  \Ob(\multc M)$, $k \geq 0$ of ``$k$-morphisms'' from $(a_1,\ldots,a_k)$
  to $b$,
\item a right action of the symmetric group $\Sigma_k$ on
  the class of all $k$-morphisms such that
  $\sigma^*$ maps the set $\multc M_k(a_1,\ldots,a_k;b)$ to the set
  $\multc M_k(a_{\sigma(1)},\ldots,a_{\sigma(k)};b)$,
\item an ``identity'' map $1_a \in M_1(a;a)$ for all $a \in \Ob(\multc
  M)$, and
\item a ``composition'' map
\begin{eqnarray*}
\multc M_n(b_1,\ldots,b_n;c) \times \multc M_{k_1}(a_{11},\ldots,a_{1 k_1};
b_1) \times \cdots\\
%\times \multc M_{k_n}(a_{n 1}, \ldots, a_{n k_n}; b_n) 
\to \multc M_{k_1 + \cdots + k_n}(a_{11}, \ldots, a_{n k_n}; c)
\end{eqnarray*}
which is associative, unital, and respects the symmetric group action.
\end{itemize}
\end{defin}

We will not make precise these last properties; they are essentially
the same as the definitions for an operad.  A map between
multicategories that preserves the appropriate structure will be
referred to as a multifunctor.

\begin{exam}
\label{exam:promote}
Any symmetric monoidal category $(\multc C, \square)$ is a multicategory, with
\[
\multc C_k(a_1,\ldots,a_k;b) = \multc C(a_1 \square \cdots \square a_k,b).
\]
For example, the categories of $\Gamma$-spaces or symmetric spectra
under $\smsh{}$ are multicategories.

There is a (lax) symmetric monoidal functor $\mb U$ from
$\Gamma$-spaces to symmetric spectra~\cite{mmss}.  This is a
multifunctor from the multicategory of $\Gamma$-spaces to
the multicategory of symmetric spectra.
\end{exam}

\begin{rmk}
The smash product of $\Gamma$-spaces of simplicial sets is
defined using left Kan extension.  As a result, we can equivalently
define a multicategory structure on $\Gamma$-spaces without
reference to the smash product by declaring the set of $k$-morphisms
from $(M_1,\ldots,M_k)$ to $N$ to be the set of collections of maps
\[
M_1(Y_1) \smsh{} \cdots \smsh{} M_k(Y_k) \to N(Y_1 \smsh{} \cdots
\smsh{} Y_k),
\]
natural in $Y_1, \ldots, Y_k$.
\end{rmk}

We will now define multicategories $\multc M$, $\multc A$, and $\multc
P$, enriched over topological spaces, as parameter multicategories for
$E_\infty$-modules, algebras, and pairings.  Let ${\cal E}(n)$ be the
space of linear isometric embeddings of ${\cal U}^{\otimes n}$ in
${\cal U}$. Together the ${\cal E}(n)$ form an $E_\infty$-operad.

\begin{defin}
The multicategory $\multc M$ has objects $R$ and $M$, such that
$\multc M_k(B_1,\ldots, B_k;C)$ is equal to ${\cal E}(k)$ in the
following cases:
\begin{itemize}
\item $B_j = C = R$ for all $j$, or
\item $B_i = C = M$ for some $i$, and $B_j = R$ for all $j \neq i$.
\end{itemize}
Otherwise, $\multc M_k(B_1, \ldots, B_k; C) = \emptyset$.  Composition
is given by the composition in ${\cal E}$.
\end{defin}

\begin{defin}
The multicategory $\multc A$ has objects $R$ and $A$, such that
$\multc A_k(B_1,\ldots, B_k;C)$ is equal to ${\cal E}(k)$ in the
following cases:
\begin{itemize}
\item $B_j = C = R$ for all $j$, or
\item $C = A$.
\end{itemize}
Otherwise, $\multc A_k(B_1,\ldots, B_k; C) = \emptyset$.  Composition
is given by composition in ${\cal E}$.
\end{defin}

\begin{defin}
The multicategory $\multc P$ has objects $R$, $M$, and $N$, and $P$,
such that $\multc P_k(B_1,\ldots,B_k;C)$ is equal to ${\cal E}(k)$
in the following cases:
\begin{itemize}
\item $B_j = C = R$ for all $j$,
\item $B_i = C = M$ for some $i$, and $B_j = R$ for all $i \neq j$,
\item $B_i = C = N$ for some $i$, and $B_j = R$ for all $i \neq j$,
\item $B_i = C = P$ for some $i$, and $B_j = R$ for all $i \neq j$, or
\item $B_i = N$, $B_{i'} = M$, $C = P$, and $B_j = R$ for
  all $j \neq i, i'$. 
\end{itemize}
Otherwise, $\multc P_k(B_1,\ldots, B_k; C) = \emptyset$.  Composition
is given by composition in ${\cal E}$.
\end{defin}

There are multifunctors $j: \multc M \to \multc A$ and $k : \multc P
\to \multc A$ with $j(R) = k(r) = R$, $j(M) = k(M) = k(N) = k(P) =
A$. Similarly, there are multifunctors $i_1, i_2, h: \multc M \to
\multc P$ that send $R$ to $R$ and such that $i_1(M) = M$, $i_2(M) =
N$, and $h(M) = P$.  These induce restriction maps on multifunctors
out to other categories; for instance, $j^*$ restricts
$E_\infty$-algebras to their underlying modules.

Additionally, all three of these multicategories have a common
subcategory $\multc R$ with a single object $R$.  Write $\eta$ for the
embedding of $\multc R$ in any of these multicategories.

\begin{prop}
\label{prop:filteredmult}
Associated to a group ${\rm G}$, there is a collection of
multifunctors as follows:
\begin{itemize}
\item $r : \multc R \to \Gamma$-spaces,
\item $m : \multc M \to \Gamma$-spaces,
\item $m_n, K_n, T_n, S_n : \multc M \to \Gamma$-spaces for $n \in \mb N$,
\item $K : \multc A \to \Gamma$-spaces, and
\item $P_{n,m} : \multc P \to \Gamma$-spaces for $n,m \in \mb N$.
\end{itemize}
These multifunctors are continuous with respect to the enrichment in
spaces, and satisfy the following properties:
\begin{itemize}
\item $r(R) = ku$,
\item $\eta^*(F) = r$ for any multifunctor $F$ on the above list,
\item $m(M) = ku$, $m_n(M) = {\cal L}(n)_+ \smsh{} ku$, where ${\cal
  L}(n)$ is the space of isometries $\mb C^n \otimes {\cal U} \to
  {\cal U}$,
\item $K(A) = {\cal K}{\rm G}$, $K_n(M) = {\cal K}{\rm G}_n$,
\item $T_n(M) = R_n \otimes_{\PU(n)} ku^{\PU(n)}$, $S_n(M)$ is the
  underlying $\Gamma$-space of $ku^{\PU(n)}$, and
\item $i_1^*(P_{n,m}) = K_n$, $i_2^*(P_{n,m}) = K_m$, $h^*(P_{n,m}) =
  K_{nm}$.
\end{itemize}
There are weak equivalences $m \leftarrow m_n \to S_n$ for each $n$,
and natural transformations $K_1 \to K_2 \to \ldots \to j^* K$ that
realize the inclusions of ${\cal K}{\rm G}_n$ into ${\cal K}{\rm G}$.  Similarly,
there are natural transformations $P_{n,m} \to P_{n',m'}$ for $n \leq
n', m \leq m'$ and $P_{n,m} \to k^*K$, which all commute, and applying
$i_1^*$, $j_1^*$, or $h^*$ yields a family of natural transformations
that realize the above inclusions.

Additionally, there are natural transformations $K_n \to T_n$ that
realize the quotient map ${\cal K}{\rm G}_n \to R_n \otimes_{\PU(n)} ku^{\PU(n)}$.
The group $\PU(n)$ acts continuously on $S_n$, and there are
$\PU(n)$-equivariant maps from $R_n$ to the set of natural
transformations ${\rm Nat}(S_n,T_n)$.  The adjoints of these maps
realize the assembly maps $R_n \smsh{\PU(n)} ku^{\PU(n)} \to R_n
\otimes_{\PU(n)} ku^{\PU(n)}$.

A map ${\rm G} \to {\rm G}'$ induces a natural transformation of
multifunctors in the opposite direction.
\end{prop}

\vskip 0.7pc

\begin{pf}
Write ${\cal K}{\rm G}^{\cal V}$ for the $\Gamma$-space ${\cal K}{\rm G}$
indexed on the universe ${\cal V}$.  For groups ${\rm G}_1, \ldots,
{\rm G}_k$, there is a well-defined exterior tensor product of
representations:
\[
{\cal K}({\rm G}_1)^{{\cal U}_1}(Z_1) \smsh{} \cdots \smsh{}
{\cal K}({\rm G}_k)^{{\cal U}_k} (Z_k) \to {\cal K}({\rm G}_1
\times \cdots \times {\rm G}_k)^{\otimes_i{\cal
    U}_i}(Z_1 \smsh{} \cdots \smsh{} Z_k).
\]
This tensor product is coherently commutative and associative with
respect to the underlying coherently commutative and associative
tensor product on inner product spaces.  It is natural in the $Z_i$,
and comes from a natural map of $\Gamma$-spaces
\[
{\cal K}({\rm G}_1)^{{\cal U}_1} \smsh{} \cdots \smsh{}
{\cal K}({\rm G}_k)^{{\cal U}_k} \to {\cal K}({\rm G}_1
\times \cdots \times {\rm G}_k)^{\otimes_i{\cal
    U}_i}.
\]
The tensor product induces coherent pairings
\[
{\cal K}({\rm G}_1)^{{\cal U}_1}_{n_1} \smsh{} \cdots \smsh{}
{\cal K}({\rm G}_k)^{{\cal U}_k}_{n_k} \to {\cal K}({\rm G}_1
\times \cdots \times {\rm G}_k)^{\otimes_i{\cal
    U}_i}_{n_1 \cdots n_k}.
\]

Now restrict to the case when ${\cal U}_i = {\cal U}$ for all $i$.
Post-composition with linear isometric embeddings ${\cal
U}^{\otimes k} \to {\cal U}$ then gives maps of $\Gamma$-spaces
\[
{\cal E}(k)_+ \smsh{} {\cal K}({\rm G}_1)_{n_1} \smsh{}
\cdots \smsh{} {\cal K}({\rm G}_k)_{n_k} \to {\cal K}({\rm G}_1 \times
\cdots \times {\rm G}_k)_{n_1 \cdots n_k}.
\]
If all ${\rm G}_i$ are equal to ${\rm G}$ or the trivial group, we
can pull back along the diagonal map to get a map
\[
{\cal E}(k)_+ \smsh{} B_1 \smsh{} \cdots \smsh{} B_k \to C,
\]
where the $B_i$ and $C$ are all of the form ${\cal K}{\rm G}$, ${\cal
K}{\rm G}_n$, or $ku$.  These maps have continuous adjoints that define the
multifunctors $r$, $m$, $K$, $K_n$, and $P_{n,m}$.  The multifunctors
$m_n$ are formed by smashing $m$ with the spaces ${\cal L}(n)_+$;
projection from ${\cal L}(n)_+ \to S^0$ gives the weak equivalence
$m_n \to m$.

Similarly, the underlying $\Gamma$-space of $ku^{\PU(n)}$ admits an
exterior tensor product.  A point of $ku^{\PU(n)}(\PU(n)_+ \smsh{} Z)$
is a map 
\[
f: Z \to \coprod V^{\cal U}(nd)/I \otimes U(d),
\]
where $V^{\cal U}(nd)$ is the Stiefel manifold of $nd$-frames in
${\cal U}$, such that $f(z) \perp f(z')$ if $z \neq z'$.  The tensor
product of frames induces a map
\[
\left[V^{\cal U}(nd)/I \otimes U(d)\right] \smsh{} \left[V^{\cal
    V}(k)/U(k)\right] \to V^{{\cal U} \otimes {\cal V}}(ndk)/I \otimes
U(dk).
\]
This product is coherently commutative and associative with
respect to the tensor product on universes.  Just as with ${\cal K}{\rm G}$,
post-composition with linear isometries gives maps
\[
{\cal E}(k)_+ \smsh{} ku \smsh{} \cdots \smsh{} ku^{\PU(n)} \smsh{} ku
\smsh{} \cdots \smsh{} ku \to ku^{\PU(n)}, 
\]
giving the desired multifunctor $S_n$.  These pairings are clearly
$\PU(n)$-equivariant, hence $\PU(n)$ acts on the functor $S_n$.

Tensoring with $\mb C^n$ and precomposing with an isometry $\mb C^n
\otimes {\cal U} \to {\cal U}$ gives a weak equivalence
\[
{\cal L}(n)_+ \smsh{} V^{\cal U}(d)/\U(d) \to V^{\cal U}(nd)/\U(nd).
\]
This commutes with the exterior tensor product, and so defines a
natural weak equivalence of multifunctors $m_n \to S_n$.

To define the multifunctor $T_n$, one can either go through the same
argument with $R_n \otimes_{\PU(n)} ku^{\PU(n)}$ or recall that it is
isomorphic to the $\Gamma$-space $F_n$.  $F_n$ is a quotient of ${\cal
K}{\rm G}_n$ by an equivalence relation that is respected by tensor product
with trivial representations.  The natural map ${\cal K}{\rm G}_n \to F_n$
yields a natural transformation of multifunctors.

The only remaining issue is to check that the assembly map is a map of
$E_\infty$-$ku$-modules.  For this, it suffices to note that the
assembly map commutes with the exterior tensor pairing
\[
(ku^{\PU(n)})^{\cal U}(\PU(n)_+ \smsh{} Y) \smsh{} ku^{\cal V}(Z)
\to ku^{\cal U \otimes V}(\PU(n)_+ \smsh{} Y \smsh{} Z)
\]
for $Y, Z \in \Gamma^o$, ${\cal U}$ and ${\cal V}$ universes.
\end{pf}

We now prove the following rigidification result.  Recall from
Example~\ref{exam:promote} that $\mb U$ is the functor which takes a
$\Gamma$-space to the associated (topological) symmetric spectrum.

\begin{prop}
\label{prop:rigid}
There is a commutative ring object in symmetric spectra $ku^r$ with module
spectra $ku^{\PU(n),r}$, and contravariant functors ${\cal K}(-)^r_n$,
${\cal K}(-)^r$ from finitely generated discrete groups to connective
$ku^r$-module symmetric spectra with the following properties:
\begin{itemize}
\item there are isomorphisms in the stable homotopy category of
  symmetric spectra $F^r \cong \Sing \mb U(F)$, where $F$ is one of
  $ku$, $ku^{\PU(n)}$, ${\cal K}(-)$, or ${\cal K}(-)_n$,
\item as a $ku^r$-module, $ku^{\PU(n),r}$ is weakly equivalent to
  $ku^r$,
\item there are $ku^r$-module maps ${\cal K}{\rm G}^r_1 \to {\cal K}{\rm G}^r_2 \to
  \cdots \to {\cal K}{\rm G}^r_\infty$, and ${\cal K}{\rm G}^r$ is weakly
  equivalent to the homotopy colimit,
\item there are strictly commutative and associative $ku^r$-module
  pairings ${\cal K}{\rm G}^r_n \smsh{ku} {\cal K}{\rm G}^r_m \to {\cal K}{\rm G}^r_{nm}$
  that commute with the above maps,
\item for any $n$, $ku^{\PU(n),r}$ is acted on by $\Sing \PU(n)$, and
  the homotopy cofiber of the map ${\cal K}{\rm G}^r_{n-1} \to {\cal
  K}{\rm G}^r_n$ is a $ku^r$-module equivalent to the derived smash product
  $\Sing R_n \smsh{\Sing \PU(n)} ku^{\PU(n),r}$.
\end{itemize}
All of the above are natural in ${\rm G}$.
\end{prop}

\begin{pf}
The multifunctors constructed in Proposition~\ref{prop:filteredmult}
are multifunctors of categories enriched in topological spaces.  Given
such a multifunctor $g: \multc C \to \Gamma$-spaces, the composite
functor $\mb U \circ g$ takes values in topological symmetric
spectra.  Similarly, the singular complex functor $\Sing$ is a
Quillen equivalence that is lax symmetric monoidal with respect to
the smash product $\smsh{}$, so there is a simplicial multifunctor
$\Sing(\mb U \circ g)$ from $\Sing \multc C$ to symmetric spectra.
Additionally, there is a simplicial multifunctor $\pi$ from $\Sing
\multc C$ to the constant simplicial multicategory $\pi_0 \multc C$.
$\pi$ is a weak equivalence.

In~\cite{rma}, a simplicial closed model structure is constructed on
the category of simplicial multifunctors from a simplicial
multicategory $\multc D$ to ${\cal S}$, the category of symmetric
spectra.  In particular, Theorem~1.4 of \cite{rma} proves that if $f :
\multc D \to \multc D'$ is a simplicial multifunctor, the restriction
map $f^* : {\cal S}^{\multc D'} \to {\cal S}^{\multc D}$ has a left
adjoint $f_*$, and if $f$ is a weak equivalence then this adjoint pair
is a Quillen equivalence.

Let $\mb L\pi_*$ denote the total left derived functor of $\pi_*$,
which consists of cofibrant resolution followed by $\pi_*$.
We obtain a rigidified symmetric spectrum $\Rig(g) =  \mb L \pi_*
\Sing(\mb U \circ g)$ such that $\pi^* \Rig(g)$ is isomorphic to
$\Sing(\mb U \circ g)$ in the homotopy category of multifunctors
$\pi_0(\multc C) \to {\cal S}$.

Additionally, in our case the multicategory $\multc C$ accepts a
``unit'' map $\eta$ from the multicategory $\multc R$ such that
$\eta^* g = r$.  Let $\tilde g \to \Sing(\mb U \circ g)$ be a
cofibrant replacement (an acyclic fibration where $\tilde g$ is a
cofibrant object) and similarly $\tilde r \to \Sing(\mb U \circ r)$.
The map $\eta^*$ preserves weak equivalences and fibrations, so the
map $\eta^* \tilde g \to \eta^* \Sing(\mb U \circ g)$ is an acyclic
fibration.  The isomorphism $r \to \eta^*g$ therefore lifts to a weak
equivalence $\tilde r \to \eta^* \tilde g$.

Composing the weak equivalences
\[
\tilde r \to \eta^* \tilde g \to \eta^* \pi^* \pi_* \tilde g \to \pi^*
\eta^* \pi_* \tilde g,
\]
we get an adjoint weak equivalence $\pi_* \tilde r \to \eta^* \pi_*
\tilde g$.  Both objects are strictly commutative ring symmetric
spectra weakly equivalent to the object $\Sing(\mb U \circ r)$.  In
other words, there is a natural weak equivalence $\eta': \Rig(r) \to \eta^*
\Rig(g)$.

Define $ku^r = \Rig(r)(R)$, ${\cal K}{\rm G}^r_n = \Rig(K_n)(M)$, ${\cal K}{\rm G}^r
= \Rig(K)(A)$, and $ku^{\PU(n),r} = \Rig(S_n)(M)$.  The unit maps
$\eta'$ make all of these objects and maps between them maps of
$ku^r$-modules.  The existence of algebra structures and  pairings of
modules are restatements of the multicategory structures on these
modules.

The weak equivalence of $ku^r$ with $ku^{\PU(n),r}$ follows from the
weak equivalences of multifunctors $m \leftarrow m_n \to S_n$.

The weak equivalence between ${\cal K}{\rm G}^r$ and the homotopy colimit of
the sequence ${\cal K}{\rm G}^r_n$ follows from Proposition~\ref{prop:hocolim},
as we have ${\cal K}{\rm G}^r_n \simeq \Sing \mb U{\cal K}{\rm G}_n$
and ${\cal K}{\rm G} \simeq \Sing \mb U{\cal K}{\rm G} \simeq \Sing
\mb U(\hocolim {\cal K}{\rm G}_n)$.

The homotopy cofiber sequence
\[
\Sing \mb U {\cal K}{\rm G}_{n-1} \to \Sing \mb U {\cal K}{\rm G}_{n} \to \Sing
\mb U F_n
\]
is weakly equivalent to the sequence of maps
\[
{\cal K}{\rm G}^r_{n-1} \to {\cal K}{\rm G}^r_n \to \Rig(T_n)(M)
\]
in the category of $ku^r$-module symmetric spectra; we now prove that
this last space is weakly equivalent to an equivariant smash product.

The $ku^r$-module $R_n \smsh{\PU(n)}
S_n(M)$ is weakly equivalent to the geometric realization of the bar
construction $B(R_n, \PU(n)_+, S_n(M))$ because $R_n$ is a free
$\PU(n)$-CW complex.  The structure map $X \smsh{} \mb UN \to \mb U(X
\smsh{} N)$ for a space $X$ and a $\Gamma$-space $N$ is an
isomorphism, and similarly $\Sing(X \smsh{} E) \to \Sing(X) \smsh{}
\Sing(E)$ is a weak equivalence for a $CW$-complex $X$ and a
topological symmetric spectrum $E$.  Therefore, both of these functors
preserve this bar construction up to weak equivalence.  The result is
that there is a weak equivalence of symmetric spectra
\begin{equation}
\label{eq:temp1}
\Sing \mb U\left(R_n \smsh{\PU(n)} S_n(M)\right) \to
\Sing(R_n) \smsh{\Sing \PU(n)} \Sing \mb U S_n(M),
\end{equation}
where this latter smash product is taken in the derived sense.

The continuous $\PU(n)$-equivariant map $R_n \to {\rm Nat}(S_n,T_n)$
induces a $\Sing \PU(n)$-equivariant map $\Sing R_n \to
F(\Rig(S_n)(M), \Rig(T_n)(M))$.  There is an adjoint map
\[
\Sing R_n \smsh{\Sing \PU(n)} \Rig(S_n)(M) \to \Rig(T_n)(M),
\]
where the smash product is taken in the derived category.  This is
weakly equivalent to the map
\[
\Sing R_n \smsh{\Sing \PU(n)} \Sing \mb U S_n(M) \to \Sing \mb U T_n(M).
\]
Using the weak equivalence of Equation \ref{eq:temp1}, this is weakly
equivalent to the map
\[
\Sing \mb U\left(R_n \smsh{\PU(n)} S_n(M)\right) \to
\Sing \mb U T_n(M).
\]
The map $R_n \smsh{\PU(n)} S_n(M) \to T_n(M)$ is a weak equivalence by
Corollary~\ref{cor:weaksmash}.  Therefore, the map $\Sing R_n
\smsh{\Sing \PU(n)} ku^{\PU(n),r} \to \Rig(T_n)(M)$ is a weak
equivalence of $ku^r$-modules.
\end{pf}

\section{The exact couple for ${\cal K}{\rm G}$}
\label{sec:exactcouple}

There is the following chain of weak equivalences of symmetric
spectra.  Here the smash products are taken in the derived category to
assure associativity.
\begin{eqnarray*}
\eilm{\mb Z} \smsh{ku^r} \left(\Sing R_n \smsh{\Sing \PU(n)}
ku^{\PU(n),r} \right) &\simeq& \Sing R_n \smsh{\Sing \PU(n)} \eilm{\mb
Z}\\ 
&\simeq& \eilm{\mb Z} \smsh{} (R_n / \PU(n)).
\end{eqnarray*}

Define $\QIrr({\rm G},n) = R_n / \PU(n)$.  $\QIrr({\rm G},n)$ is the
quotient space of isomorphism classes of representations ${\rm G}$ of
dimension $n$ modulo decomposable representations.  (The notation is
to avoid confusion with the standard notation for the {\em subspace\/}
of isomorphism classes of irreducible representations.)

Proposition~\ref{prop:rigid} identifies the following homotopy cofiber
sequences.  The homotopy colimit of the top row is weakly equivalent
to ${\cal K}{\rm G}^r$.

\xym{
{*} \ar[r] & {\cal K}{\rm G}^r_1 \ar[d] \ar[r] & {\cal K}{\rm G}^r_2 \ar[d] \ar[r]
& \ldots\\
& \Sing R_1 \smsh{} ku^r
& \Sing R_2 \smsh{\Sing \PU(2)} ku^{\PU(2),r}
}

Using the weak equivalence ${\cal K}{\rm G}^r \simeq \Sing(\mb U {\cal
K}{\rm G})$, the following spectral sequence results.

\begin{thm}
\label{thm:kspecseq}
There exists a convergent right-half-plane spectral sequence of the
form
\[
E_1^{p,q} = ku^{\PU(n)}_{q-p+1}(R_{p-1}) \Rightarrow
\pi_{p+q}({\cal K}{\rm G}).
\]
\end{thm}

\begin{rmk}
This uses a Serre indexing convention, so that $d_r$ maps $E^{p,q}_r$
to $E^{p-r,q+r-1}_r$.  We should remark that by the homotopy group
$\pi_n$ of a symmetric spectrum $X$, we mean the set $[S^n,X]$ in the
stable homotopy category, or equivalently the stable homotopy groups
of an $\Omega$-spectrum weakly equivalent to $X$.  The symmetric
spectra ${\cal K}{\rm G}$ and ${\cal K}{\rm G}_n$ are almost $\Omega$-spectra, and
so in this case the notion coincides with the classical notion of
stable homotopy groups.
\end{rmk}

Smashing the previous diagram over $ku^r$ with $\eilm{\mb Z}$ (with
smash product taken in the derived category of $ku^r$-modules) yields
the following.  The homotopy colimit of the top row is weakly
equivalent to $\eilm{\mb Z} \smsh{ku^r} {\cal K}{\rm G}^r$.

\xym{
{*} \ar[r] & \eilm{\mb Z} \smsh{ku^r} {\cal K}{\rm G}^r_1 \ar[d] \ar[r] & \eilm{\mb
  Z} \smsh{ku^r} {\cal K}{\rm G}^r_2 \ar[d] \ar[r] & \eilm{\mb Z}
\smsh{ku^r} {\cal K}{\rm G}^r_3 \ar[d] \ldots\\
& \eilm{\mb Z} \smsh{} {\QIrr({\rm G},1)} & \eilm{\mb Z} \smsh{}
\QIrr({\rm G},2) & \eilm{\mb Z} \smsh{} \QIrr({\rm G},3)
}

Again, the diagram results in a spectral sequence.

\begin{thm}
\label{thm:homspecseq}
There exists a convergent right-half-plane spectral sequence of the
form
\[
E_1^{p,q} = H_{q-p+1}(\QIrr({\rm G},p-1)) \Rightarrow \pi_{p+q}(\eilm{\mb
  Z} \smsh{ku^r} {\cal K}{\rm G}).
\]
\end{thm}

\begin{exam}
When ${\rm G}$ is finite or nilpotent, the cofiber sequences are all
split.  When ${\rm G}$ is finite, this is clear.  When ${\rm G}$ is
nilpotent, results of~\cite{lubotzkymagid} show that the space of
irreducible representations of dimension $n$ is closed in $\Hom({\rm
G},\U(n))$, which provides the desired splitting $R_n \smsh{\PU(n)}
ku^{\PU(n)} \to {\cal K}{\rm G}_n$.

As a result, for these groups we have a weak equivalence
\[
{\cal K}{\rm G} \simeq \bigvee \left(R_n \smsh{\PU(n)} ku^{\PU(n)}\right).
\]
In this case, the spectral sequence of Theorem~\ref{thm:homspecseq}
degenerates at the $E_1$ page.  For example, consider the integer
Heisenberg group of $3 \times 3$ strict upper triangular matrices with
integer entries.  The structure of the space of all irreducible
representations appears in~\cite{nunleymagid}.

The $E_1 = E_\infty$ page of the spectral sequence of
Theorem~\ref{thm:homspecseq} looks, in part, as follows:

\[
\xy 0;/r1.5pc/:
(0,0);(7,0)**\dir{-};
(0,-4);(0,4.8)**\dir{-};
(0.5,0.5),*{\mb Z},
c+(1,-1),*{\mb Z},
c+(1,-1),*{\mb Z},
c+(1,-1),*{\mb Z},
(0.5,1.5),*{\mb Z^2},
c+(1,-1),*{\mb Z^2},
c+(1,-1),*{\mb Z^2},
c+(1,-1),*{\mb Z^2},
c+(1,-1),*{\mb Z^2},
(0.5,2.5),*{\mb Z},
c+(1,-1),*{\mb Z},
c+(1,-1),*{\mb Z},
c+(1,-1),*{\mb Z},
c+(1,-1),*{\mb Z},
c+(1,-1),*{\mb Z},
(0.5,3.5),*{0},c+(1,-1),*{0},c+(1,-1),*{0},c+(1,-1),*{0},c+(1,-1),*{0},c+(1,-1),*{0},
(1.5,3.5),*{0},c+(1,-1),*{0},c+(1,-1),*{0},c+(1,-1),*{0},c+(1,-1),*{0},
(2.5,3.5),*{0},c+(1,-1),*{0},c+(1,-1),*{0},c+(1,-1),*{0},
(3.5,3.5),*{0},c+(1,-1),*{0},c+(1,-1),*{0},
(4.5,3.5),*{0},c+(1,-1),*{0},
(5.5,3.5),*{0},
(3.5,4.5),*{\vdots},
(3.5,-3.5),*{\vdots},
(6.5,0.5),*{\cdots}
\endxy
\]
\end{exam}

\begin{exam}
Suppose ${\rm G}$ is free on $k$ generators.  As mentioned in the
introduction, there is a weak equivalence ${\cal
K}{\rm G} \simeq ku \vee (\vee^k \Sigma ku)$.

When ${\rm G}$ is free on two generators, explicit computations with
the spectral sequence of Theorem~\ref{thm:homspecseq} give the
following picture of the $E_1$ page.

\[
\xy 0;/r1.5pc/:
(0,0);(5,0)**\dir{-};
(0,-3);(0,7)**\dir{-};
(2.5,-1.5),*{0},c+(1,0),*{0},
(1.5,-0.5),*{0},c+(1,0),*{0},c+(1,0),*{?},
(0.5,0.5),*{\mb Z},c+(1,0),*{0},c+(1,0),*{?},c+(1,0),*{?},
(0.5,1.5),*{\mb Z^2},c+(1,0),*{0},c+(1,0),*{?},c+(1,0),*{?},
(0.5,2.5),*{\mb Z},c+(1,0),*{\mb Z},c+(1,0),*{?},c+(1,0),*{?}, 
(0.5,3.5),*{0},c+(1,0),*{\mb Z^2},c+(1,0),*{?},c+(1,0),*{?},
(0.5,4.5),*{0},c+(1,0),*{\mb Z},c+(1,0),*{?},c+(1,0),*{?},
(0.5,5.5),*{0},c+(1,0),*{0},c+(1,0),*{?},c+(1,0),*{?},
(1.2,2.5);(0.8,2.5)**\dir{-}*\dir{>},
(3.5,6.5),*{\vdots},
(3.5,-2.5),*{\vdots},
(4.5,0.5),*{\cdots}
\endxy
\]

The differential $d_1: E_1^{1,2} \to E_1^{0,2}$ is an isomorphism.
The terms $E_1^{p,q}$ are zero on the set $\{p > 0, p + q < 2\}$, and
also on the set $\{q > p^2 + p + 2\}$.  The terms where $q = p^2 + p +
2$ are all isomorphic to $\mb Z$.

This spectral sequence converges to $\mb Z$ in dimension 0, $\mb Z^2$
in dimension 1, and $0$ in all other dimensions.  The classes in 
$E_1^{0,0}$ and $E_1^{0,1}$ are precisely those classes that survive
to the $E_\infty$ term.
\end{exam}

\section{Proof of the product formula}
\label{sec:product}

In this section we will prove Theorem~\ref{thm:product}, the product
formula for deformation $K$-theory spectra.

The proof requires the following lemmas.
\begin{lem}
\label{lem:connectivity}
A map $M' \to M$ of connective $ku^r$-module spectra is a weak
equivalence if and only if the map $\eilm{\mb Z} \smsh{ku^r} M' \to
\eilm{\mb Z} \smsh{ku^r} M$ is a weak equivalence.  (Smash products
are taken in the derived category.)
\end{lem}

\begin{pf}
By taking cofibers, it suffices to prove the equivalent statement that
a connective $ku^r$-module spectrum $M''$ is weakly contractible if and
only if $\eilm{\mb Z} \smsh{ku^r} M'' \simeq *$.

However, smashing $M''$ with the homotopy cofiber sequence $\Sigma^2 ku^r
\to ku^r \to \eilm{\mb Z}$ of $ku^r$-module spectra shows that
$\eilm{\mb Z} \smsh{ku^r} M'' \simeq *$ if and only if the Bott map
$\beta : \Sigma^2 M'' \to M''$ is a weak equivalence.  This would
imply that the homotopy groups of $M''$ are periodic; $M''$ is
connective, so the result follows.
\end{pf}

\begin{lem}
\label{lem:prodrep}
Irreducible unitary representations of ${\rm G} \times {\rm H}$ are
precisely of the form $V \otimes W$ for $V, W$ irreducible unitary
representations of ${\rm G}$ and ${\rm H}$ respectively.
\end{lem}

\begin{pf}
That the tensor product of irreducible group representations is
irreducible, and conversely irreducible representations of ${\rm G}
\times {\rm H}$ are tensor products, is well known.  (See, for
example, \cite{vinberg} Section 4.4, Theorem 6 and the comment
afterwards.)

Suppose $V$ and $W$ are nontrivial irreducible representations of
${\rm G}$ and ${\rm H}$ respectively.  Clearly invariant inner
products on $V$ and $W$ induce one on $V \otimes W$.  Conversely, if
$W$ is nonzero $V$ appears as a sub-${\rm G}$-representation of $V
\otimes W$, and hence an invariant inner product on $V \otimes W$
induces one on the subspace $V$.  Therefore, $V \otimes W$ admits a
unitary structure if and only if both $V$ and $W$ do.
\end{pf}

\begin{rmk}
Lemma~\ref{lem:prodrep} fails when we consider representations of
the group ${{\rm G}\times{\rm H}}$ in other groups such as
orthogonal groups and symmetric groups.
\end{rmk}

\begin{pf}[of Theorem~\ref{thm:product}]
The proof consists of constructing a filtration of the spectrum ${\cal
K}{\rm G}^r \smsh{ku^r} {\cal K}{\rm H}^r$ that agrees with the existing
filtration on ${\cal K}({\rm G} \times {\rm H})^r$.

We apply the results of Proposition~\ref{prop:rigid} to get a map of
$ku^r$-algebras
\[
{\cal K}{\rm G}^r \smsh{ku^r} {\cal K}{\rm H}^r \to {\cal K}({\rm G} \times {\rm H})^r
\smsh{ku^r} {\cal K}({\rm G} \times {\rm H})^r \to {\cal K}({\rm G}
\times {\rm H})^r.
\]
Similarly, whenever $p \cdot q \leq n$ there is a corresponding map of
$ku^r$-modules
\[
{\cal K}{\rm G}^r_p \smsh{ku^r} {\cal K}{\rm H}^r_q \to {\cal
  K}({\rm G} \times {\rm H})^r_n.
\]
This diagram is natural in $p$, $q$, and $n$.

Let $\Gamma$ denote a cofibrant replacement functor for
$ku^r$-modules.  If we define new $ku^r$-module spectra $M_n =
\hocolim_{p \cdot q \leq n} \Gamma{\cal K}{\rm G}^r_p \smsh{ku^r} \Gamma{\cal
K}{\rm H}^r_q$, then there are induced $ku$-module maps
\[
f_n : M_n \to {\cal K}({\rm G} \times {\rm H})^r_n.
\]

The maps $(\hocolim {\cal K}{\rm G}^r_p) \to {\cal K}{\rm G}^r$ and $(\hocolim
{\cal K}{\rm H}^r_q) \to {\cal K}{\rm H}^r$ are weak equivalences, so there is a weak
equivalence
\[
\hocolim M_n \simeq \hocolim_{p,q} \Gamma{\cal K}{\rm G}^r_p \smsh{ku^r}
\Gamma{\cal K}{\rm H}^r_q \simeq \Gamma{\cal K}{\rm G}^r \smsh{ku} \Gamma{\cal K}{\rm H}^r.
\]
Therefore, it suffices to show $M_n \to {\cal K}({\rm G} \times {\rm
H})^r_n$ is a weak equivalence for all $n$.  We have an induced map of
homotopy cofiber sequences:

\xym{
M_{n-1} \ar[r] \ar[d]_{f_{n-1}} & 
M_n \ar[r] \ar[d]_{f_n} & 
M_n/M_{n-1} \ar[d]_{g_n}
\\
{\cal K}({\rm G} \times {\rm H})_{n-1}^r \ar[r]&
{\cal K}({\rm G} \times {\rm H})_n^r \ar[r]&
R_n({\rm G} \times {\rm H}) \smsh{\PU(n)} ku^{\PU(n),r}.
}

To prove the theorem it suffices to show that the map $g_n$ is a weak
equivalence for all $n$.

The spectra $M_n/M_{n-1}$ and $R_n({\rm G} \times {\rm H})
\smsh{\PU(n)} ku^{\PU(n),r}$ are connective $ku^r$-module spectra.
Applying Lemma~\ref{lem:connectivity}, it suffices to prove that the
map $\eilm{\mb Z} \smsh{ku^r} g_n$ is a weak equivalence.

Because the map $M_{n-1} \to M_n$ is a map from the homotopy colimit
of a subdiagram into the full diagram, we can explicitly compute the
homotopy cofiber of this map.  The homotopy cofiber is weakly
equivalent to the wedge
\[
\bigvee_{p \cdot q = n} \left(\Gamma{\cal K}{\rm G}^r_p / \Gamma{\cal
    K}{\rm G}^r_{p-1}\right) \smsh{ku^r} \left(\Gamma{\cal K}{\rm H}^r_q /
  \Gamma{\cal K}{\rm H}^r_{q-1}\right).
\]
To see this, one notes that we can take the homotopy colimit by
replacing the diagram by a homotopy equivalent diagram made up of
cofibrations, and then take the take the ordinary colimit.  For any
such diagram $\{F_{p,q}\}$, we clearly have
\[
\left(\bigcup_{p\cdot q \leq n} F_{p,q}\right) \Big/
\left(\bigcup_{p\cdot q < n} F_{p,q}\right) = \bigvee_{p \cdot q = n}
F_{p,q} / \left(F_{p-1,q} \cup_{F_{p-1,q-1}} F_{p,q-1}\right).
\]
The ``pushout product axiom'' \cite{ss} shows that in the case where
$F_{p,q} \simeq A_p \smsh{} B_q$ and the maps $A_{p-1} \to A_p$ and
$B_{p-1} \to B_p$ are cofibrations, we have
\[
F_{p,q} / \left(F_{p-1,q} \cup_{F_{p-1,q-1}} F_{p,q-1}\right) \simeq 
A_p/A_{p-1} \smsh{ku^r} B_q/B_{q-1}.
\]

The spectra ${\cal K}{\rm G}^r_p / {\cal K}{\rm G}^r_{p-1}$, and the
corresponding 
spectra for ${\rm H}$, are weakly equivalent to those that were
identified as equivariant smash product spectra in
Corollary~\ref{cor:fibseq} and Proposition~\ref{prop:baseid}.
Smashing over $ku^r$ with $\eilm{\mb Z}$ gives us the following
identity.
\begin{eqnarray*}
\eilm{\mb Z} \smsh{ku^r} M_n / M_{n-1} &\simeq& \bigvee_{p \cdot q = n}
\left(\eilm{\mb Z} \smsh{} \QIrr({\rm G},p)\right) \smsh{\eilm{\mb Z}}
\left(\eilm{\mb Z} \smsh{} \QIrr({\rm H},q)\right) \\
&\simeq& \eilm{\mb Z} \smsh{} \left(\bigvee_{p \cdot q = n} \QIrr({\rm
    G},p) \smsh{} \QIrr({\rm H},q)\right).
\end{eqnarray*}

The map $\eilm{\mb Z} \smsh{ku} g_n$ can be identified with the map
\[
\eilm{\mb Z} \smsh{} \left(\vee_{p \cdot q=n} \QIrr({\rm G},p) \smsh{}
\QIrr({\rm H},q)\right) \to \eilm{\mb Z} \smsh{} \QIrr({\rm G}
\times {\rm H},n)
\]
that is induced by the tensor product of representations.  The tensor
product map $\otimes: \vee_{p \cdot q = n} \QIrr({\rm G},p) \smsh{}
\QIrr({\rm H},q) \to \QIrr({\rm G} \times {\rm H},n)$ is a
continuous map between compact Hausdorff spaces.  It is bijective by
Lemma~\ref{lem:prodrep}.  Therefore, it is a homeomorphism.
\end{pf}

\begin{acknowledgements}
The author would like to extend thanks Gunnar Carlsson, Chris Douglas,
Bjorn Dundas, Haynes Miller, and Daniel Ramras for many helpful
conversations, and to Mike Hill for his comments.
\end{acknowledgements}

\nocite{*}
\bibliography{lawson_pf}

\end{article}
\end{document}